\title{\bf
Convergence of Nonlinear Observers on $\reals^n$
with a Riemannian Metric (Part I)
}
\author{Ricardo G. Sanfelice\thanks{R. G. Sanfelice is with
the Department of Aerospace and Mechanical Engineering, University of Arizona
1130 N. Mountain Ave, AZ 85721,
       Email: {\tt\small sricardo@u.arizona.edu}}
and Laurent Praly\thanks{L. Praly is with CAS, ParisTech, Ecole des Mines, 35 rue Saint Honor\'{e}, 77305, 
        Fontainebleau, France
       Email: {\tt\small Laurent.Praly@ensmp.fr}
}}

\documentclass[journal,onecolumn,draft]{IEEEtran} 

\usepackage{hyperref}
\usepackage{amsmath}
\usepackage{amssymb}
\usepackage{latexsym}
\usepackage[usenames]{color} 
\usepackage{psfrag}
\usepackage{graphics}
\usepackage{graphicx}
\usepackage{setspace}
\usepackage{balance}

 \definecolor{mygreen}{rgb}{0,.7,.9}

\definecolor{MyMagenta}{rgb}{0.8,0,0.8} 

\newcommand{\non}{\nonumber}

\def\XX{X}
\def\sphere{{\mathbb{S}}}
\newcommand{\reals}{\mathbb{R}}
\newcommand{\A}{{\cal{A}}}
\newcommand{\bigbrace}[1]{\left\{#1\right\}}
\newcommand{\defset}[2]{\bigbrace{#1\ : \ #2 }}
\newcommand{\x}{x}
\newcommand{\xdot}{\dot{\x}}
\newcommand{\f}{f}
\newcommand{\h}{h}
\newcommand{\y}{y}
\newcommand{\ox}{\chi}
\newcommand{\oxdot}{\dot{\ox}}
\newcommand{\F}{F}
\renewcommand{\H}{H}

\newcommand{\xhat}{{\hat{x}}}
\newcommand{\xhatdot}{\dot{\xhat}}

\newcommand{\OurProblem}{($\star$)}

\newcommand{\matt}[1]{\begin{bmatrix}#1\end{bmatrix}} 
\newcommand{\RR}{{\mathbb R}}

\catcode`\@=11
\def\downparenfill{$\m@th\braceld\leaders\vrule\hfill\bracerd$}
\def\overparen#1{\mathop{\vbox{\ialign{##\crcr\crcr \noalign{\kern0.4ex}
\downparenfill\crcr\noalign{\kern0.4ex\nointerlineskip}
$\hfil\displaystyle{#1}\hfil$\crcr}}}\limits}
\catcode`\@=12

\def\dotoverparen#1{\mathop{\vbox{\ialign{##\crcr\hfill$\cdot $\hfill\crcr 
\noalign{\kern-3.7ex}
\downparenfill\crcr\noalign{\kern0.4ex\nointerlineskip}
$\hfil\displaystyle{#1}\hfil$\crcr}}}\limits}

\newcommand{\realsgeq}{[0,+\infty)}

\newcommand{\plower}{\underline{p}}
\newcommand{\pupper}{\overline{p}}
\newcommand{\qlower}{q}

\def\Did{{\mathfrak{D}^+}}

\def\dcurve{\omega}

\usepackage{ifthen}

\newtheorem{theorem}{Theorem}[section] 
   
\newtheorem{lemma}[theorem]{Lemma}   
\newtheorem{example}[theorem]{Example}    
\newtheorem{proposition}[theorem]{Proposition}    
\newtheorem{remark}[theorem]{Remark}   
\newtheorem{definition}[theorem]{Definition}

\begin{document}
\maketitle

\vspace{-0.5in}
\begin{abstract}
We study how convergence of an observer whose state lives in a copy of the given system's space can be established using a Riemannian metric.  We show that the existence of an observer guaranteeing the property that a Riemannian distance between system and observer solutions is nonincreasing implies that the Lie derivative of the Riemannian metric along the system vector field is conditionally negative. Moreover, we establish that the existence of this metric is related to the observability of the system's linearization along its solutions. Moreover, if the observer has an infinite gain margin then
the level sets of the output function are geodesically convex. Conversely, we establish that, if a complete Riemannian metric has a Lie derivative along the system vector field that is conditionally negative and is such that the output function has a monotonicity property, then there exists an observer with an infinite gain margin.
\end{abstract}

\section{Introduction}
\label{sec:Introduction}
For a nonlinear system of the form
\begin{equation}
\label{eqn:Plant1}
\xdot  \;=\;  \f(\x)\  ,\quad 
\y  \;=\;  \h(\x)
\end{equation}
with $\x \in \reals^n$ being the system's state and $\y\in\reals^m$ the 
measured system's output,
we study the problem of obtaining an estimate $\xhat$ 
of the state $\x$ by means of the dynamical system, called {\em observer},
\begin{equation}
\label{eqn:Observer1}
\oxdot \;=\;  \F(\ox,\y)\  ,\quad 
\xhat \;=\;  \H(\ox,\y)
\end{equation}
with $\ox \in \reals^p$ being the observer's state and
$\xhat \in \reals ^n$ the observer's output, 
used as the system's state estimate.
We focus on the case where the
state $\ox$ of the observer evolves in
a copy of the space of the system's state $\x$, i.e.,
they both belong to $\reals^n$,
with,
moreover,
an
output function $H$ such that
$\xhat = \ox$.
We consider the following observer design problem:
\begin{itemize}
\item[\OurProblem] 
Given functions
$\f$ and $\h$,
design a function $\F$ such that
for
the system
\begin{eqnarray}\label{eqn:Interconnection}
\begin{array}{lll}
\xdot  &=&  \f(\x)\  ,\quad 
\xhatdot \ = \  \F(\xhat,\h(\x)),
\end{array}
\end{eqnarray}
the
zero estimation error
set
\begin{equation}\label{eqn:GASset}
\A\;=\; 
\defset{(\x, \xhat)\in\reals^n\times\reals^{n}}{ x = \xhat}
\end{equation}
is
globally asymptotically stable
(see the text below \eqref{eqn:KLBoundError}).
\end{itemize}

Many contributions from different viewpoints have been made
to address problem \OurProblem.
While a summary of the very rich literature on the topic is out of 
the scope of this paper, it is important to point out
the interest of exploiting a possible contraction
property of the flow generated by the observer.
Study of contracting flows has a very long 
history and has been proposed independently by several authors; see, e.g.,
\cite{Lewis.71,Hartmann.64,Demidovich.61-62,Nemeth.98,Lohmiller.Slotine.98.Automatica} 
(see \cite{Jouffroy.05.CDC} for a historical discussion).
In the context of observers, 
Riemannian metrics have been used in \cite{Aghannan.Rouchon.03,Bonnabel.07,Bonnabel.10},
for instance,
with the objective of guaranteeing that the Riemannian distance
between the system and observer solutions decreases to zero. 
In these papers,
the authors consider systems whose dynamics follow from
a principle of least action involving a Riemannian metric, 
such as Lagrangian systems with a Lagrangian that is quadratic in the generalized 
velocities.
The observer design therein exploits some properties of this 
metric and local convergence is established via some ad-hoc modification of 
this metric
or choice of coordinates.

This paper advocates that,
since the observability of the system linearized along each of its solutions may vary 
significantly from one solution to another, 
the native Euclidean geometry of the state space
may not be appropriate to study 
convergence properties of an observer.
Instead of insisting in using a Riemannian metric associated to the system's dynamics, 
we propose to study Riemannian metrics incorporating information on the
system's dynamics and observability.
In Section \ref{sec1b}, we show that
if for a given Riemannian metric
an observer whose state $\ox$ lives  in
a copy of
the given system's state space and
makes the Riemannian distance along
system and observer solutions nonincreasing
then, necessarily, the Lie derivative of the metric along the system 
solutions
satisfies an inequality involving 
the output function.
Section \ref{sec1c} shows that
if the same conditions hold and the observer has an infinite gain 
margin then, necessarily,
the level sets of the output 
function are geodesically convex. 
In Section \ref{sec1d}
we establish that if
a Riemannian metric with a Lie derivative satisfying the inequality 
mentioned above is, in some coordinates,
uniformly bounded away from zero and upper bounded
then the system's linearization along each of its solution must be 
detectable.
With
the insight provided by 
these necessary conditions, 
Section \ref{sec3} proposes 
a set of sufficient conditions guaranteeing the existence of an 
observer whose flow leads to a decreasing Riemannian distance between
system's state and estimated state.

For the sake of simplicity, we assume throughout the paper that the functions are differentiable sufficiently many times.
Moreover, we work under restrictions that can be further relaxed, such as time independence of the
right-hand sides and  
forward 
completeness of the systems\footnote{A system is said to be {\em forward complete} if each of its solutions exists on $[0,+\infty)$.}.

This paper is devoted to analysis.
In a companion paper,
we focus on observer design, namely, 
on the construction of a Riemannian metric
satisfying the desired inequality on its Lie derivative and making 
the level sets of the output function possibly totally geodesic.

\begin{example}[Motivational example]
\label{ex:1}
We illustrate our 
results in the
following
academic system
\begin{equation}
\label{LP64}
\dot x_1 = x_2\,  \sqrt{1+x_1^2},\quad 
\dot x_2\;=\; - \frac{x_1}{\sqrt{1+x_1^2}}\,  x_2^2,\quad 
y\;=\; x_1
\  .
\end{equation}
For this system
 (\ref{LP64}), by following \cite{Krener.Isidori.83}, we get the observer
\begin{equation}
\label{LP70}
\begin{array}{l}
\displaystyle
\dot {\bar {\hat{x}}}_1\;=\; \bar {\hat{x}}_2 - (\bar {\hat{x}}_1-y),\quad 
\dot {\bar {\hat{x}}}_2\;=\; - (\bar {\hat{x}}_1-y),\\ \displaystyle
 \hat x_1 \;=\; \bar {\hat{x}}_1,\quad  \hat x_2\;=\; \frac{\bar {\hat{x}}_2}{\sqrt{1+y^2}}
\  .
\end{array}
\end{equation}
This observer is in the form (\ref{eqn:Observer1}), but cannot be 
written in the form of (\ref{eqn:Interconnection}) with the $(\hat x_1,\hat x_2)$ 
coordinates since this would involve $x_2$. Nevertheless,
with the
Lyapunov function
\begin{equation}\label{eqn:Vexample}
\begin{array}{ll}\displaystyle
V(\xhat,x)\;=\; & (\xhat_1-x_1)^2\;-\; (\xhat_1-x_1)\,  (\xhat_2-x_2)\,  \sqrt{1+x_1^2}\\
& \displaystyle
\;+\; 
(\xhat_2-x_2)^2\,  (1+x_1^2)
\end{array}
\end{equation}
we obtain for the system-observer interconnection
 (\ref{LP64})-(\ref{LP70})
$$
\dotoverparen{V(\xhat,x)}
\;=\; -V(\xhat,x)
\  .
$$
Since $V$ satisfies,
for all $(x,\xhat) \in \reals^2\times\reals^{2}$,
\begin{eqnarray*} \displaystyle
\frac{(\xhat_1-x_1)^2+(\xhat_2-x_2)^2}{2}
 \leq V(\xhat,x) & \\ & \hspace{-1.2in} \leq
 \displaystyle
\frac{3}{2}
\left[(\xhat_1-x_1)^2+(\xhat_2-x_2)^2 \right]
\left(1 + x_1^2 \right),
\end{eqnarray*}
this implies that,
for all $t \geq 0$ and all $(x,\xhat) \in \reals^2\times\reals^{2}$,
\begin{equation}\label{eqn:KLBoundError}
| X(\x,t) - \hat X((\hat x,\x),t) |^2 
\leq
3
\exp(-t) 
(1+ x_1^2) 
| \x - \xhat|^2 
\   ,
\end{equation}
where
$(X(\x,t),\hat X((\hat x,\x),t))$
is the solution issued from points $(\x,\xhat)$
for the system-observer interconnection (\ref{LP64})-(\ref{LP70}).
This establishes that the set $\A$ is
globally  asymptotically stable 
(nonuniformly in $x$ but 
uniformly in  $x - \hat x$).

As it will be shown in Section~\ref{sub:Riemannian},
the key point here is that
$V$ is the square of a Riemannian distance between $\xhat$ and $x$
that is associated to an $x$-dependent Riemannian
metric. Moreover, as justified in Section~\ref{sec1b}, 
no matter what the observer is,
it is impossible to find
a standard quadratic form 
expressed in the given coordinates
(i.e., a Riemannian distance
associated with a constant Riemannian metric)
that is nonincreasing along solutions.
This is a motivation for the analysis of observers
using $x$-dependent Riemannian metrics.
\hfill$\Box$
\end{example}

\section{Necessary conditions for having a Riemannian distance 
between system and observer solutions to decrease.}
\label{sec1}

\subsection{Riemannian Distance}
\label{sub:Riemannian}

As discussed in Section~\ref{sec:Introduction}, 
the notions
of nonexpanding/contracting flow and geodesically monotone vector 
fields
are suitable
for studying
asymptotic stability of the zero error set $\A$
in (\ref{eqn:GASset}). 
We start by recalling some basic facts on Riemannian 
distance.

Let 
$P:\reals^n\to\reals^{n\times n}$
be a 
$C^3$
symmetric covariant two-tensor (see, e.g., \cite[Page 17]{Sakai.96}).
If $x$ and $\bar x$ are two sets 
of coordinates related by
$
\bar x \;=\; \phi(x)
$
with 
$\phi$ 
being
a diffeomorphism,
then $P$ expressed in $x$ coordinates as $P(x)$
and in $\bar x$ coordinates as ${\bar P}({\bar x})$
are related by (see, e.g., \cite[Example II.2]{Sakai.96})
\begin{equation}
\label{LP69}
P(x)\;=\; \frac{\partial \phi}{\partial x}(x)^\top \bar P(\bar x) \,  \frac{\partial 
\phi}{\partial x}(x)
\  .
\end{equation}
 If
$P$ takes positive definite values
then 
the length of a $C^1$ path
$\gamma $ between points $x_1$ and $x_2$ 
is defined as
\begin{equation}\label{eqn:Ldefn}
\left.\vrule height 1em depth 0.5em width 0pt
L(\gamma )\right|_{s_1}^{s_2}\;=\; \int_{s_1}^{s_2}
\sqrt{
\frac{d\gamma }{ds}(s) ^{\top}P(\gamma (s)) \frac{d\gamma }{ds}(s)
} \,  ds,
\end{equation}
where
$$
\gamma (s_1)\;=\; x_1\quad ,\qquad \gamma (s_2)\;=\; x_2
\  .
$$
With such a definition, $P$ is also called a Riemannian metric.
The Riemannian distance $d(x_1,x_2)$ is 
the minimum of $\left.\vrule height 1em depth 0.5em width 0pt
L(\gamma )\right|_{s_1}^{s_2}$ among all possible piecewise $C^1$ 
paths $\gamma $ between $x_1$ and $x_2$.
To relate the Riemannian distance with geodesics,
we invoke
the Hopf-Rinow Theorem 
(see, e.g., \cite[Theorem II.1.1]{Sakai.96}),
which asserts the following:
if every geodesic
can be maximally extended to $\reals$ then the minimum of 
$\left.\vrule height 1em depth 0.5em width 0pt
L(\gamma )\right|_{s_1}^{s_2}$
is actually given 
by the length of a 
(maybe nonunique) geodesic, which is called a {\em minimal geodesic};
for more details, see, e.g., \cite{Boothby.75} and \cite{DoCarmo.92}.
In the appendix we show that, in our context,
this maximal extension property holds on 
$\reals^n$
if
there exist globally defined coordinates 
in which $P$ satisfies
\begin{equation}
\label{LP45}
0\: <\:  P(x)
\quad \forall x\in \reals ^n\  ,\quad 
\lim_{r \to \infty }r^2 \plower(r) \:=\: +\infty,
\end{equation}
where, for any positive real number $r$,
$$
\plower(r)\;=\; \min_{x:|x| \leq r}
\lambda _{\min}\left(P(x)\right)
\  ,
$$
with $\lambda _{\min}\left(P(x)\right)$ denoting
the minimum eigenvalue of $P(x)$.
In this case, the Riemannian metric given by $P$ is said to be complete
and, denoting by $\gamma ^*$ a minimal
(normalized\footnote{
A normalized geodesic $\gamma^*$ satisfies $\  \frac{d\gamma^* }{ds}(s)^\top 
P(\gamma^* (s))\,  \frac{d\gamma^* }{ds}(s)\   = 1$ for all $s$ in its domain 
of definition. In the following, the adjective ``normalized'' is 
omitted.
})
geodesic between $x=\gamma ^*(0)$
and $\xhat = \gamma ^*(\hat s) $,
with $\hat s \geq 0$,
the Riemannian distance $d(\xhat,\x)$ is
\begin{equation}\label{eqn:Rdistance}
d(\xhat,\x)\;=\; \left.
\vrule height 1em depth 0.5em width 0pt
L(\gamma^*)\right|_{0}^{\hat s}\;=\;
\hat s
\  .
\end{equation}

\begin{example}
As an illustration,
consider the symmetric covariant two-tensor expressed in $x$ coordinates
as
$
P(x) = 
 \matt{
1 -\frac{x_1x_2}{\sqrt{1 + x_1^2}} + \frac{x_1^2 x_2^2}{1+x_1^2}
&
-\frac{\sqrt{1 + x_1^2}}{2} + x_1 x_2
\\[0.7em]
-\frac{\sqrt{1 + x_1^2}}{2} + x_1 x_2
&
1+x_1^2
}$.
Since
condition \eqref{LP45} holds with $\underline{p}(r) = \frac{1}{2}$
for all $r > 0$,
it is a complete Riemannian metric.
Moreover, using \eqref{LP69}, it is easy to check that in the coordinates
$
\bar x = \phi(x) = \matt{x_1 \\ x_2 \sqrt{1 + x_1^2}},\ \mbox{ its expression is }
\bar P(\bar x) = 
\matt{1 &\displaystyle -\frac{1}{2} \\
\displaystyle -\frac{1}{2} & 1}.
$
Since $\bar P(\bar x)$ is constant, 
any minimal geodesic $\bar \gamma^*$
takes the form 
$\bar \gamma^*(s) = \bar x + s \bar v
$
with $\bar v \in \reals^2$ satisfying $\bar v^\top \bar P(\bar x) \bar v = 1$.
Then, 
a minimal geodesic in $x$ coordinates
is given by
$\gamma^*(s) = \phi^{-1}(\bar x + s \bar v)$.
Accordingly, the Riemannian distance between $\xhat$ and $x$ is
\begin{eqnarray} \label{eqn:VisDforExample}
\begin{array}{lll}
\int_{0}^{\hat s}
\sqrt{
\frac{d\gamma^* }{ds}(s) ^{\top}P(\gamma^* (s)) \frac{d\gamma^* }{ds}(s)
} \,  ds
&=&  d(\hat { x},x) \ = \ d(\bar {\hat x},\bar x)
\\ \non
& &\hspace{-1.2in} = \int_{0}^{\hat s}
\sqrt{
\frac{d\bar \gamma^* }{ds}(s) ^{\top}\bar P(\bar \gamma^* (s)) \frac{d\bar \gamma^* }{ds}(s)
} \,  ds \\
& &\hspace{-1.2in} = 
\sqrt{
\left(\bar {\hat x} - \bar x)^\top \bar P (\bar x) (\bar {\hat x} - \bar x\right)}
\\ \non
& &\hspace{-1.2in} = 
\sqrt{
( \phi (\hat x) - \phi (x))^\top \bar P (\bar x) (\phi (\hat x) - \phi (x))} \\
& &\hspace{-1.2in} =  \sqrt{V(\xhat,x)},
 \end{array}
\end{eqnarray}
where $V$ is given in \eqref{eqn:Vexample}
and 
$\bar {\hat x} = \phi(\hat x)$.
\hfill$\Box$
\end{example}
\par\vspace{1em}

Having a Riemannian distance, we say that a system
$
\dot x =  f(x),
$
with solutions $X(x,t)$, generates a 
nonexpanding 
(respectively, contracting) flow if, for any pair $(x_1,x_2)$ in 
$\reals^n\times\reals^n$, the function
$t\mapsto d(X(x_1,t),X(x_2,t))$
is nonincreasing (respectively, 
strictly decreasing);
see, e.g., \cite{Isac.Nemeth.08}.
Also, the vector field $f$ is said 
to be
geodesically
monotonic (respectively, strictly monotonic) if we have
\begin{equation}
\label{LP51}
\mathcal{L}_f P(x) \; \leq \; 0
\qquad (\textrm{respectively}\  \; <\; 0)
\qquad \forall x\in \reals^n
\  ,
\end{equation}
where $\mathcal{L}_f P$ 
is the Lie derivative 
of the symmetric covariant two-tensor 
$P$, whose expression in $x$ coordinates is
\begin{equation}\label{LP23}
\begin{array}{lll}
\displaystyle v^\top \mathcal{L}_f P(x)\,  v &  & \\
& & \hspace{-0.9in}= \displaystyle
\lim_{r\to 0} 
\left[\frac{[(I+r\frac{\partial f}{\partial x}(x))v]^\top 
P(x+rf(x)) [(I+r\frac{\partial f}{\partial x}(x))v]}{r}\right.\\[0.7em]
& & \hspace{-0.5in}\displaystyle \left.- \frac{v^\top P(x) v}{r}\right]\\
& & \hspace{-0.8in} = \displaystyle
\frac{\partial}{ \partial \x}\left(
\vrule height 1em depth 0.5em width 0pt
v^{\top} P(\x) \,  v\right)\,  \f(x) \;+\; 2\,  
v^{\top} P(\x)\left(\frac{\partial \f}{\partial \x}(\x) \,  v\right)
\end{array}
\end{equation}
for all $v\in \reals^n$;
see \cite[Exercise V.2.8]{Boothby.75}, \cite[Page 17]{Sakai.96}, or \cite{Lebedev.Cloud.05}.
We have
the following result 
(see, for instance, \cite{Isac.Nemeth.08} or \cite{Aghannan.Rouchon.03} for 
a proof).
\begin{lemma}
{\it
A geodesically monotonic (respectively, strictly monotonic) vector field generates a
nonexpanding (respectively, contracting) flow.
}
\end{lemma}
If inequality (\ref{LP51})
holds for the observer vector field $F$ then 
$t\mapsto d(\hat X((\xhat _1,x),t),\hat X((\xhat _2,x),t))$
is (respectively, strictly) decreasing;
however,
this property is more than 
what is needed for the zero estimation error set $\A$
to be (respectively, asymptotically) stable.
Actually, it is sufficient 
to have an observer giving rise
to a (respectively, strictly) decreasing function
$t\mapsto d(\hat X((\xhat,x),t),X(x,t))$ for all pairs $(\xhat,x)$ in $\reals^n\times\reals^n$.
That is, we
do not insist on having a Riemannian distance between any two arbitrary observer solutions
to decrease, but only to have a decreasing Riemannian distance between any observer
solution and its corresponding system solution (which is a particular observer solution).

\subsection{Necessity of geodesic monotonicity in the directions tangent to the level sets of the output function}
\label{sec1b}
Since the Riemannian distance 
between $\xhat$ and $x$
is locally Lipschitz,
its upper right-hand Dini derivative is given by
\begin{equation}
\label{LP57}%
\Did d(\xhat,\x)\;=\; 
\limsup_{t\to 0_+} \frac{d(\hat X((\hat x,\x),t),X(x,t))-d(\hat x,x)}{t}
\end{equation}
for each $(\xhat,x) \in \reals^n \times \reals^n$.
It is nonpositive
when the function
$t\mapsto d(\hat X((\xhat,x),t),X(x,t))$
is nonincreasing.

\begin{theorem}
\label{thm3}
\textit{
Assume there exists a complete $C^3$ Riemannian metric $P$ such that,
for each $(\xhat,x) \in \reals^n \times \reals^n$,
\begin{equation}\label{eqn:Ddot1}
\Did d(\xhat,\x) \leq 0\ 
\end{equation}
holds along any solution of 
\eqref{eqn:Interconnection},
then 
\begin{equation}\label{eqn:LfPDnonPositive}
\begin{array}{ll}
v^\top \mathcal{L}_f P(x) v
\;\leq \; 0\,   &\qquad 
\forall (x,v) \in \reals^n\times\reals^n\\
& \qquad \mbox{ such that }\ \ \displaystyle \frac{\partial h}{\partial x}(x)v = 0\ .
\end{array}
\end{equation}
Furthermore, if 
there exists a function
$\dcurve:\reals^n \times \reals^n \to \realsgeq$
such that $(\xhat,x)\mapsto d(\xhat,x) \dcurve(\xhat,x)$
is a $C^2$ function on a neighborhood ${\mathcal N}_\A$ of $\A$ 
with the property that, for some $\varepsilon > 0$,
\begin{equation}\label{eqn:D2dcurveBound}
\frac{\partial^2 (d\, \dcurve) }{\partial \xhat^2}(x,x) \geq \varepsilon P(x) \qquad \forall x \in 
\reals^n
\end{equation}
and, for each $(\xhat,x) \in {\mathcal N}_\A$,
\begin{equation}\label{eqn:Ddot2}
\Did d(\xhat,\x) \leq - \dcurve(\xhat,x)
\end{equation}
holds along any solution of 
\eqref{eqn:Interconnection}, then there exists a continuous function 
$\rho :\reals^n\to \RR$ satisfying
\begin{equation}
\label{LP71}
\mathcal{L}_f P(x)\;\leq \; \rho (x)\,  
\frac{\partial \h}{\partial \x}(\x)^{\top} \frac{\partial \h}{\partial 
\x}(\x) - \frac{\varepsilon}{2}P(x)
\qquad
 \forall x \in \reals^n.
\end{equation}
}
\end{theorem}
\begin{IEEEproof}
To simplify the notation,
let $V:\reals^n\times\reals^n\to \realsgeq$
be the function defined as
the square of the Riemannian distance, i.e.,
$V(\xhat,x) = d(\xhat,x)^2$,
and notice that\footnote{
Since
$
\limsup (a\, b) \leq \limsup a \cdot \limsup b.
$
}
\begin{eqnarray}\label{eqn:DdotProperty}
\Did V(\xhat,\x)  =  \Did d^2(\xhat,\x) \leq 2\, d(\xhat,x)\, \Did d(\xhat,\x).
\end{eqnarray}

Pick an arbitrary point $x$ in $\reals^n$.
From \cite[Theorem 3.6]{Kaboyashi.Nomizu.96},
there exists a (normal coordinate) neighborhood 
$\mathcal{N}_x$ such that $V$  is $C^2$ 
on $\mathcal{N}_x\times \mathcal{N}_x$.
From \eqref{eqn:DdotProperty}
and \eqref{eqn:Ddot1}
(respectively, 
from \eqref{eqn:DdotProperty}
and \eqref{eqn:Ddot2},
on $({\mathcal N}_x \times {\mathcal N}_x) \cap {\mathcal N}_\A$),
we have
\begin{eqnarray*}
\Did V(\xhat,\x) \leq 0 \qquad (\, \mbox{respectively } \leq - 2\, d(\xhat,x) \dcurve(\xhat,x)\, ).
\end{eqnarray*}
Let $r_*$ be
a strictly real number 
such that, for any $v$ in $\sphere^n$, the unit sphere, and for all 
$r \in [0,r_*)$, $(\hat x+rv,x)$ are the 
coordinates of a point in $(\mathcal{N}_x \times \mathcal{N}_x) \cap \mathcal{N}_\A$.
We have
\footnote{
This follows from the fact that a first order approximation of the geodesic is
$
\gamma (s) = x + s\, v+ O_{x,v}(s^2)
$
with
$
v^\top P(x)\,  v\;=\; 1
$,
which yields
$
V(\xhat,x)=d(\xhat,x)^2\;=\; \hat s^2\;=\; (\xhat -x)^\top P(x)\,  (\xhat 
-x)\;+\; O_{x,v}(\hat s^3)$, 
where 
the subindex in $O_{\x,v}$ indicates dependence on $(\x,v)$.
}
\begin{equation}\label{eqn:Pfromd2dx2}
\frac{\partial ^2 V}{\partial \xhat^2}(x,x)=
\frac{\partial ^2 V}{\partial x^2}(x,x)\;=\; 
2
P(x)
\end{equation}
and\footnote{%
This follows from $\x=\xhat$ being a minimizer of $V$ for all $x$.
}
\begin{equation}
\label{LP1}
\renewcommand{\arraystretch}{1.5}
\begin{array}{@{}c@{}}
\displaystyle
V(\x,\x) \:=\:  0
\  ,
\displaystyle 
\frac{\partial V}{\partial x}(\x,\x)\:=\:
\frac{\partial V}{\partial \xhat}(\x,\x)\;=\; 0
\  
\\\displaystyle 
\frac{\partial ^2 V}{\partial x^2}(\x,\x) +
\frac{\partial ^2 V}{\partial \xhat\partial x}(\x,\x)
\!\!\;=\;\!\!
\frac{\partial ^2 V}{\partial \xhat^2}(\x,\x)\;+\;
\frac{\partial ^2 V}{\partial x\partial \xhat }(\x,\x)
\:=\:  0
\end{array}
\end{equation}
and, for all $r \in [0,r^*)$ and $v \in \sphere^n$,
\begin{equation}\non
\begin{array}{lll}
\Did V(x+rv,\x) &=& \displaystyle
\frac{\partial V}{\partial \x}(x+rv,x) \,  f(\x) \\
& & \;+\; \displaystyle
\frac{\partial V}{\partial \xhat}(x+rv,x)\,  F(x+rv,h(x))\\
&  \leq &  0 \\
(\  \mbox{respectively}\  
&\leq& \displaystyle\; - 2\, d(x+rv,x)\, \dcurve(x+rv,x)\  ).
\end{array}
\end{equation}
With the definition of $d$, this implies that 
$\A$ is forward invariant, i.e.,
the solutions to
\eqref{eqn:Interconnection}
with
$x=\xhat$
as initial 
condition remain in $\A$ for all $t \geq 0$. 
This implies
\begin{eqnarray}\label{eqn:FandfEquality}
\F(\x,\h(\x)) = \f(\x)
\  .
\end{eqnarray}
By differentiating this identity with respect to $\x$, we get
\begin{eqnarray}\label{eqn:dF}
\frac{\partial \F}{\partial \x}(\x,\h(\x)) 
\;+\;  \frac{\partial \F}{\partial y}(\x,h(\x)) \frac{\partial \h}{\partial \x}(\x) 
\;=\; \frac{\partial \f}{\partial \x}(\x).
\end{eqnarray}

For $r$ in $(0,r_*)$,
we obtain
\\[1em]$\null\quad\displaystyle 
\frac{1}{r^2} \left[\frac{\partial V}{\partial \x}(\x+r v,\x) + \frac{\partial V}{\partial \xhat}(\x+r v,\x)
\right]\f(\x) +$\\[0.7em]$\null\hfill\displaystyle
\frac{\displaystyle\frac{\partial V}{\partial \xhat}(\x+r v,\x)}{r} 
\frac{\F(\x + rv,\h(\x))-f(x)}{r} \leq 0
$\refstepcounter{equation}\label{eqn:Limit}\hfill$(\theequation)$
\\[0.7em]\null\hfill$\displaystyle 
\qquad (\mbox{respectively}
\leq \; -  \frac{2}{r^2} d(\x+r v,\x)\, \dcurve(\x+r v,\x)).
\quad$\\[1em]
To compute the limit for $r$ approaching $0$
note that we have the following Taylor expansion
 around 
$(\x,\x)$
\begin{eqnarray*}
V(\x+r v,\x) &=&
V(\x,\x) + r\,  
\frac{\partial V}{\partial \xhat }(\x,\x) \,  v\\
& & 
+ \frac{r^2}{2}\,  
v^{\top}
\frac{\partial^2 V}{\partial \xhat ^2}(\x,\x)\,  v + 
O_{\x,v}(r^3)\  ,
\end{eqnarray*}
\begin{eqnarray*}
\frac{\partial V}{\partial \xhat}(\x+r v,\x) &=&
\frac{\partial V}{\partial \xhat}(\x,\x) + r\,  
\frac{\partial^2 V}{\partial \xhat^2 }(\x,\x) \,  v
+ 
O_{\x,v}(r^2)\  ,
\end{eqnarray*}
\begin{eqnarray*}
\frac{F(x+rv,h(x))-f(x)}{r} & = & \frac{F(x,h(x))-f(x)}{r}\\
& & + \frac{\partial F}{\partial 
\xhat}(x,h(x))\,  v
\;+\; O_{\x,v}(r).
\end{eqnarray*}
Define $W(x) = V(x+r\, v, x)$ 
and note that
\begin{eqnarray*}
\frac{\partial W}{\partial x}(x) =
\frac{\partial V}{\partial \xhat}(\x+r v,\x)
+
\frac{\partial V}{\partial \x}(\x+r v,\x).
\end{eqnarray*}
With  \eqref{eqn:Pfromd2dx2} and \eqref{LP1}, we get
\begin{eqnarray*}
W(\x) &=&
r^2\,  
v^{\top}
P(\x)\,  v + 
O_{\x,v}(r^3)\  ,\\
\frac{1}{r}\frac{\partial V}{\partial \xhat}(\x+r v,\x) &=&
2\,  
v^\top P(\x) 
+ 
O_{\x,v}(r),
\end{eqnarray*}
and with \eqref{eqn:FandfEquality}
\begin{eqnarray*}
\frac{F(x+rv,h(x))-f(x)}{r} & = & \frac{\partial F}{\partial 
\xhat}(x,h(x))\,  v
\;+\; O_{\x,v}(r).
\end{eqnarray*}
This yields
\begin{equation}
\label{eqn:Limit2}
\begin{array}{lll}
\displaystyle
\lim_{r\to 0} \frac{1}{r^2}\left[
\frac{\partial V}{\partial \xhat}(\x+r v,\x)
+
\frac{\partial V}{\partial \x}(\x+r v,\x)
\right] f(x)& & \\
 & & \hspace{-2.5in}= \displaystyle
\lim_{r\to 0} \frac{1}{r^2} 
\frac{\partial W}{\partial \x}(x) f(x)
=
\frac{\partial \left(v^\top P\,  v\right)}{\partial 
x}(x) f(x).
\end{array}
\end{equation}
Also, with (\ref{eqn:FandfEquality}), we get
\begin{equation}
\label{eqn:Limit3}
\begin{array}{lll}
\displaystyle
\lim_{r\to 0} 
\frac{\displaystyle\frac{\partial V}{\partial \xhat}(\x+r v,\x)}{r} 
\frac{\F(\x + rv ,\h(\x))-f(x)}{r}
& & \\
& & \hspace{-1.7in} \displaystyle
 = 2
 v^\top P(x)\,  \frac{\partial F}{\partial x}(x,h(x))\,  v.
\end{array}
\end{equation}
Similarly, we can obtain
\begin{equation}
\label{LP2}
\lim_{r \to 0} \frac{2}{r^2} d(\x+r v,\x)\,\dcurve(\x+r v,\x)\;=\;
\,  v^\top \frac{\partial^2 (d\,\dcurve)}{\partial \xhat^2}(x,x)\,  v
\  .
\end{equation}

Then, combining \eqref{eqn:Limit2}, \eqref{eqn:Limit3},
and (\ref{LP2}),
we have that inequality \eqref{eqn:Limit}
gives
\begin{eqnarray*}
\frac{\partial \left(  v ^{\top}  P  v \right)}{\partial \x}(x)
\f(\x) 
+ 
2 v ^{\top} \!P(x) \frac{\partial \F}{\partial \x}(\x,\h(\x))   v 
&\leq&   0 \\
(\mbox{respectively }
\leq   -   
v^\top \frac{\partial^2 (d\,\dcurve)}{\partial \xhat^2}(x,x) v
& & \displaystyle 
\hspace{-0.2in}  \forall v \in \sphere^n 
\ ),
\end{eqnarray*}
or, equivalently, using \eqref{eqn:dF} and (\ref{LP23}),
\begin{eqnarray}\label{eqn:LfP1}
v^\top \mathcal{L}_f P(x)  v -
2v^\top P(x) \frac{\partial \F}{\partial \y}(\x,\h(\x)) \frac{\partial \h}{\partial \x}(\x)   v
&\hskip-0.5em \leq &\hskip-0.5em 0 \\ \label{eqn:LfP2}
\hspace{-0.2in}(\mbox{respectively }   \leq  
- 
v^\top \frac{\partial^2 (d\,\dcurve)}{\partial \xhat^2}(x,x) v
\ \ \forall v \in \sphere^n).
& & \displaystyle
\end{eqnarray}
It follows that \eqref{eqn:LfP1} already implies \eqref{eqn:LfPDnonPositive}.
Also, when (\ref{eqn:Ddot2}) holds, by
completing squares and using Cauchy-Schwarz inequality,
we get successively, for any
function
$\rho:\reals^n \to (0,+\infty)$ and all 
$(x,v) $ in $ \reals^n \times\sphere^n$,
\begin{eqnarray*}
2v^\top P(x) \frac{\partial \F}{\partial \y}(\x,\h(\x)) \frac{\partial \h}{\partial \x}(\x)   v 
\;
&\leq& \;  \displaystyle 
\!\!\!\rho(x) 
\left|
\frac{\partial \h}{\partial \x}(\x)   v
\right|
^2\\ 
& & \hspace{-0.7in}
+
\frac{1}{\rho(x)}
 \left|
 v^{\top} P(x) \frac{\partial F}{\partial \y}(\x,h(\x))
 \right|
 ^2
\\
&\leq& \; \displaystyle 
\!\!\!\rho(x) v^{\top} \frac{\partial \h}{\partial \x}(\x)^{\top} \frac{\partial \h}{\partial \x}(\x)   v
\\
& & \hspace{-1.3in}+
\frac{
\left|
\frac{\partial F}{\partial \y}(\x,h(\x))^\top P(x)\frac{\partial F}{\partial \y}(\x,h(\x))
\right|
}{\rho(x)} \,  v^{\top} P(x) v
\  .
\end{eqnarray*}
Equation (\ref{LP71}) follows from 
\eqref{eqn:D2dcurveBound} by picking $\rho $ as any 
continuous
function satisfying

$$
\frac{2}{\varepsilon}
\left|
\frac{\partial F}{\partial \y}(\x,h(\x))^\top P(x)\frac{\partial F}{\partial \y}(\x,h(\x))
\right|
\: \leq \: \rho(x)$$
for all $x \in \reals^n$.
\end{IEEEproof}

When compared with \eqref{LP51}, which says
 $f$ is (respectively, strictly) geodesically monotonic,
 the necessary condition
\eqref{eqn:LfPDnonPositive} (respectively, \eqref{LP71})
says only that the 
vector field $f$ is
 geodesically (respectively, strictly) monotonic
in the 
directions $v$ satisfying $\frac{\partial h}{\partial x}(x) v=0$, 
i.e., in the directions tangent to the level sets of the output 
function $h$. 

\begin{remark}
\normalfont
Theorem~\ref{thm3} can be interpreted as an extension of
\cite[Proposition 3]{Praly.01.NOLCOS.Observers}.
In this reference,
a $C^{\infty}$ function $V$ depending only on $\xhat-x$,
called a {\em state-independent error Lyapunov function},
is obtained from stability properties of $\A$.
In such a  case, the conditions in \eqref{LP1} 
yield a constant matrix $P$.
Then, Theorem~\ref{thm3} implies that, for all $x \in \reals^n$, 
$P$ is a semidefinite positive matrix that satisfies, for all $x \in \reals^n$,
$$
P\, \frac{\partial  f}{\partial x}(x)
\;+\; 
 \frac{\partial  f}{\partial x}(x)^\top P
\;\leq \; \rho (x)\,  
\frac{\partial \h}{\partial \x}(\x)^{\top} \frac{\partial \h}{\partial 
\x}(\x) 
\;-\;  \frac{\varepsilon}{2}\,  P
\  .
$$
It follows that, for all $x \in \reals^n$ and $c \in [0,\frac{\varepsilon}{4}]$, 
we have the implication
\begin{equation}\label{eqn:NOLCOSextension}
\frac{\partial h}{\partial x}(x)v = 0 \quad \Rightarrow
\quad v^{\top}P\, \frac{\partial  f}{\partial x}(x)v\;\leq \; 
\;-\;  c\, v^{\top}Pv
\  .
\end{equation}
When $c = 0$, this property corresponds to the one established in \cite[Proposition 3]{Praly.01.NOLCOS.Observers}.
It is worth pointing out that a limitation of the work in \cite{Praly.01.NOLCOS.Observers}
is that
the results are extrinsic, i.e., they depend
on the coordinates
since a quadratic form may not be 
quadratic after a nonlinear change of coordinates.
On the other hand, 
the necessary conditions in Theorem~\ref{thm3}
are intrinsic.  In fact,
let $\phi$ be a diffeomorphism on $\reals^n$ leading to the new 
coordinates
\begin{equation}\label{eqn:Diffeomorphism}
\bar x\;=\; \phi(x)\quad ,\qquad \bar {\hat {x}}\;=\; \phi(\xhat)
\  .
\end{equation}
Let $\bar h$, $\overline{d}, \overline{\dcurve},$ $\bar \rho$, $\bar f$, and $\bar P$  
be $h$, $d$, $\dcurve$, $\rho$, $f$, and $P$,
respectively, in the new coordinates. 
We have \eqref{LP69} and
\begin{eqnarray}\label{eqn:FunctionsInNewCoordinates}
\renewcommand{\arraystretch}{2}
\begin{array}{lll}\displaystyle
\bar h(\bar x)=h(x)\displaystyle
\  ,\ \
\frac{\partial h}{\partial x}(x)\ =\
\frac{\partial \bar h}{\partial \bar x}(\bar x)
\frac{\partial \phi}{\partial x}(x),\\\displaystyle
\bar f( \bar x)  =  \displaystyle\frac{\partial \phi}{\partial x}(x)\,  f(x)
\  ,\ \ \displaystyle\\
\bar d(\bar{\hat x},\bar x)=d(\hat x,x),\qquad \bar \omega(\bar{\hat x},\bar x)=\omega(\hat x,x)
\\
\displaystyle
\frac{\partial^2 (d\,\dcurve)}{\partial \xhat^2}(x,x)  = \displaystyle \frac{\partial \phi}{\partial x}(x)^\top 
\frac{\partial^2 (\overline{d}\, \overline{\dcurve})}{\partial \bar{\hat{x}}^2}(\bar x,\bar x)\frac{\partial \phi}{\partial x}(x) \, ,
\\ \non\displaystyle
\bar \rho(\bar x) = \rho(x)
\  ,
\ \
\mathcal{L}_f P(x) \ = \ \frac{\partial \phi}{\partial x}(x)^\top \mathcal{L}_{\bar f} \bar P(\bar x) \frac{\partial \phi}{\partial x}(x)
\  .
\end{array}
\end{eqnarray}
Substituting these expressions in \eqref{LP71}, we get
\begin{eqnarray*}
\frac{\partial \phi}{\partial x}(x)^\top \mathcal{L}_{\bar f} \bar P(\bar x) \frac{\partial \phi}{\partial x}(x)
\!\!\! & \leq&  \!\!\! 
\bar \rho(\bar x)\,  
\left[\frac{\partial \bar h}{\partial \bar x}(\bar x)
\frac{\partial \phi}{\partial x}(x)
\right]
^{\top}\times \\
& & \hspace{-1.1in} 
\left[\frac{\partial \bar h}{\partial \bar x}(\bar x)
\frac{\partial \phi}{\partial x}(x)
\right]
 - 
 \frac{1}{2}
\frac{\partial \phi}{\partial x}(x)^\top 
\frac{\partial^2 (\overline{d}\, \overline{\dcurve})}{\partial \bar{\hat{x}}^2}(\bar x,\bar x)\frac{\partial \phi}{\partial x}(x)
\end{eqnarray*}
and since 
$\frac{\partial \phi}{\partial x}(x)$ is invertible 
it gives
\begin{eqnarray*}
 \mathcal{L}_{\bar f} \bar P(\bar x) 
\; & \leq&  \; 
\bar \rho(\bar x)\,  
\frac{\partial \bar h}{\partial \bar x}(\bar x)
^{\top} 
\frac{\partial \bar h}{\partial \bar x}(\bar x)
 - 
 \frac{1}{2}
\frac{\partial^2 (\overline{d}\, \overline{\dcurve})}{\partial \bar{\hat{x}}^2}(\bar x,\bar x),
\end{eqnarray*}
which is  inequality (\ref{LP71}) in $\bar x$ coordinates.

Furthermore, from the definition of 
$\mathcal{L}_f P$
and with completion of squares as in the proof of Theorem~\ref{thm3},
it can be checked that 
condition \eqref{LP71} is preserved, but with a modified function $\rho$, 
after an output-dependent time scaling of the system, i.e.,
when $f$ is replaced by $\bar f(x) = \theta(h(x)) f(x)$
with $\theta$ taking strictly positive values.\hfill$\Box$
\end{remark}

The necessary conditions in Theorem~\ref{thm3} can be used to characterize 
the family of Riemannian metrics possibly leading
to a Riemannian distance that is nonincreasing (via \eqref{eqn:LfPDnonPositive}) 
or strictly decreasing (via \eqref{LP71}) along solutions.
For instance, condition \eqref{eqn:LfPDnonPositive} 
can be used to justify that, for system
\eqref{LP64},
there is no 
such a Riemannian metric that is constant.
\begin{example}[Motivational example -- continued]
\label{ex:1-re}
For the family of constant Riemannian metrics of the form
$
P =  \matt{
p & q \\
q & r
}, p, r > 0 \ , p\, r >q^2 
$
for
\eqref{LP64}, for each $v \in \reals^2$ such that
$$\frac{\partial h}{\partial x}(x)v = \matt{1 & 0} \matt{v_1 \\ v_2} =  0
\  ,
$$
we obtain
\begin{eqnarray*}
v^\top P\, \frac{\partial  f}{\partial x}(x) v
\;+\; 
v^\top  \frac{\partial  f}{\partial x}(x)^\top P v
 & & 
 \\
 & & \hspace{-1in}= \frac{2}{\sqrt{1+x_1^2}}\:    v ^\top  P\matt{ x_1 x_2 &  1+ x_1^2  \\  \displaystyle - \frac{x_2^2}{1+x_1^2} & - 2 x_1 x_2  } v
\\
  &  & \hspace{-0.8in} =
\frac{v_2^2 (2\, q\, (1+ x_1^2) -4\, r\, x_1\, x_2)
}{\sqrt{1+x_1^2}}\  ,
\end{eqnarray*}
which cannot be nonpositive for each $x$.
On the other hand, 
it can be shown that the family of Riemannian metrics
satisfying \eqref{eqn:LfPDnonPositive} 
can be described as
\begin{equation}\label{eqn:PfamilyExample}
P(x) =
\matt{1 & \frac{ x_1  x_2}{\sqrt{1+x_1^2}} \\ 0 & \sqrt{1+x_1^2} }
\matt{
\bar p({\bar x}) & \bar q({\bar x}) \\
\bar q({\bar x}) & \bar r({\bar x})
}
\matt{1 & 0 \\ \frac{ x_1  x_2}{\sqrt{1+x_1^2}} & \sqrt{1+x_1^2} }
\end{equation}
with $({\bar{x}}_1, {\bar{x}}_2 ) = (x_1,x_2\, \sqrt{1+x_1^2})$ and
$\bar r({\bar x}) = a({\bar x})^2$, 
$\bar q({\bar x}) = -b({\bar x})^2 - \frac{1}{2}\frac{\partial \bar r}{\partial {\bar x}_1}({\bar x}) {\bar x}_2$,
$\bar p({\bar x}) = c({\bar x})^2 + \frac{\bar q({\bar x})^2}{\bar r({\bar x})}$,
where $a, b, c: \reals^2 \to \reals$ are sufficiently smooth functions
with $a$ and $c$ not vanishing.
A particular choice is
$a(\bar{x}) = 1$,
$b(\bar{x}) = \frac{1}{(1+\bar x_1^2)^{\frac{1}{4}}}$,
and $c(\bar{x})^2 = 1 + \left(
\frac{\bar x_2}{1+\bar x_1^2} + \frac{\bar x_1}{\sqrt{1+\bar x_1^2}}
\right)^2$,
which leads to 
\begin{equation}\label{eqn:PexampleInx}
P(x)\;=\; \matt{
2+x_2^2& x_1x_2-1 \\ x_1x_2 -1& 1+x_1^2
}.
\end{equation}\hfill\null $\Box$
\end{example}

\subsection{
Necessity of geodesic convexity of the level sets of the output function
}
\label{sec1c}

In Theorem~\ref{thm3}, we studied 
the implications of the
existence of an observer making
$t\mapsto d(\hat X((\xhat,x),t),X(x,t))$ nonincreasing, 
in particular, when $\xhat$ converges to $x$
(in the proof, $(x+rv,x)$ approaches $(x,x)$). 
Now we study 
the implications of the
existence of such an observer 
for the case
when $\xhat$ is far 
away from $x$. To this end,
for each $s$ in $[0,\hat s]$, let
$t\mapsto\Gamma (s,t)$ be a $C^1$ function satisfying
\begin{eqnarray*}
&
\displaystyle \frac{\partial X}{\partial t}(x,t)\;=\; f(X(x,t))
\  ,\quad  
X(x,0)\;=\; x
\  ,
\\
&\displaystyle \frac{\partial \hat X}{\partial t}(\xhat,t)\;=\; F(\hat X(\xhat,t),h(X(x,t)))
\  ,\quad  
\hat X(\xhat,0)\;=\; \hat x
\  ,
\\
&
\displaystyle \frac{\partial \Gamma }{\partial t}(s,t)\;=\; F(\Gamma (s,t),h(X(x,t)))
\  ,\quad  
\Gamma (s,0)\;=\; \gamma ^*(s)
\  ,
\end{eqnarray*}
with $\gamma ^*$ a minimal geodesic between $x$ and $\hat x$.
Then, we have
$
\hat X((\hat x,\x),t)\;=\; \Gamma (\hat s,t)
$
and
hence, at time $t$, $s\mapsto \Gamma(s,t)$ is a path between $X(x,t)$ and $\hat X((\hat x,\x),t)$. 
Also,
we have
$$
d(\xhat,x)\;=\; d(\Gamma (\hat s,0),\Gamma (0,0))\; = \; \left.
\vrule height 1em depth 0.5em width 0pt
L(\Gamma (\mbox{\LARGE .}\, ,0))\right|_{0}^{\hat s}
\  .
$$
Also, we know from the first
order variation formula (see, for instance, 
\cite[Theorem 6.14]{Spivak.79} or
\cite[Theorem 5.7]{Isac.Nemeth.08})
that we have
\begin{eqnarray*}
\left.\frac{d}{dt} \left.
\vrule height 1em depth 0.5em width 0pt
L(\Gamma (\mbox{\LARGE .}\, ,t))\right|_{0}^{\hat s}\right|_{t=0}
&  &  \\
& & \hspace{-0.8in}=
\left.\frac{d}{dt}
\int_{0}^{\hat s}
\sqrt{
\frac{\partial \Gamma  }{\partial s}(s,t) ^\top P(\Gamma (s,t)) 
\frac{\partial \Gamma }{\partial s}(s,t)
} \;   ds\right|_{t=0}
\\
& & \hspace{-0.8in}=
\frac{d\gamma ^*}{ds}(\hat s)^\top P(\gamma ^*(\hat s))  \,  \F(\gamma ^*(\hat s),\y) \\
& & \hspace{-0.6in}
- \frac{d\gamma ^*}{ds}(0)^\top P(\gamma ^*(0))\,  \F(\gamma ^*(0),\y).
\end{eqnarray*}
On the other hand, in general, 
for each $t$ in the domain of definition,
we have only
$$
d(\hat X(\xhat ,t),X(x,t))\;=\; 
d(\Gamma (\hat s,t),\Gamma (0,t))\; \leq \; \left.
\vrule height 1em depth 0.5em width 0pt
L(\Gamma (\mbox{\LARGE .}\, ,t))\right|_{0}^{\hat s}
\  .
$$
Then,
the upper right-hand Dini derivative
of the distance between $\xhat$ and $\x$
 in
\eqref{LP57}
satisfies
\begin{eqnarray}
\nonumber
\Did d(\xhat,\x) &\leq &
\left.\frac{d}{dt} \left.
\vrule height 1em depth 0.5em width 0pt
L(\Gamma (\mbox{\LARGE .}\, ,t))\right|_{0}^{\hat s}\right|_{t=0}
\\ \non
&\leq &
\frac{d\gamma ^*}{ds}(\hat s)^\top P(\gamma ^*(\hat s))  \,  \F(\gamma ^*(\hat s),\y)\\
& & \label{LP7}
- \frac{d\gamma ^*}{ds}(0)^\top P(\gamma ^*(0))\,  f(\gamma ^*(0))
\  .
\end{eqnarray}
Even though \eqref{LP7} is an inequality condition,
we proceed as if it were 
an equality.
In such a case, if
the observer makes the distance $d(\xhat,x)$ nonincreasing
along solutions then necessarily the right-hand side of 
(\ref{LP7}) has to be nonpositive.
To get a better understanding of what this means,
consider the case when\footnote{
For a given $x\in\reals^n$,
this condition holds for every minimal geodesic $\gamma^*$ such that
$\frac{d\gamma ^*}{ds}(0)$ belongs to the
closed half space 
$\{w \in\reals^n : w^{\top}P(x)f(x) \leq 0\}$.
}
\begin{equation}
\label{LP44}
- \frac{d\gamma ^*}{ds}(0)^\top P(\gamma ^*(0))\,  f(\gamma ^*(0)) \; \geq 
 \; 0
\  .
\end{equation}
Then,
for the right-hand side of (\ref{LP7}) to be
nonpositive, with $\xhat=\gamma ^*(\hat s)$, we must have
\begin{equation}
\label{LP43}
\frac{d\gamma ^*}{ds}(\hat s)^\top P(\xhat)  \F(\xhat,\y)
\; \leq \; 0
\  .
\end{equation}
At this point, it is important to note that $\frac{d\gamma ^*}{ds}(\hat 
s)$ is the direction in which the 
state estimate $\xhat$ ``sees'' the system state $\x$ along
a minimal geodesic. Such a direction is unknown to the observer.
The only known information is
that, 
for given $y$,
$x$ belongs to the following 
$y$-level set\footnote{
By 
$y$-level set of $h$ we mean the intersection, for each $i = 1,2,\ldots, m$, of the sets
$
 \{x\in\reals^n : h_i(x) =y_i\}.
$
} of the 
output function:
$$
\mathfrak{H}(y)\;=\; \{x: h(x) =y\}
\  .
$$
Hence,
(\ref{LP43}) implies the following property:
given $\xhat$ and $y$, the level set of the 
output function $\mathfrak{H}(y)$ is ``seen'' from $\xhat$
along a minimal geodesic,
within a cone whose aperture is less than $\pi $.
As stated in Lemma \ref{lem2} below, this
property implies that $\mathfrak{H}(y)$
is geodesically convex;
see \cite[Definition 6.1.1]{Rapcsak.97}
and \cite[Section 9.4]{Hicks.65}.

\begin{definition}[geodesic convexity]
\label{def:GeodesicConvexity}
{\it
A subset $S$ of $\,\reals^n$ is said to be geodesically convex
if, for any pair  of points $(x_1,x_2) \in
S\times S
$, there exists a minimal geodesic 
$\gamma ^*$
between $x_1=\gamma ^*(s_1)$ and $x_2=\gamma ^*(s_2)$
satisfying
$$
\gamma ^*(s)\in S
\qquad \forall s\in [s_1,s_2].
$$
}
\end{definition}

\begin{lemma}
\label{lem2}
{\it
Let $P:\reals^n\to\reals^{n\times n}$ be a complete Riemannian metric. 
Assume $S$ is
a subset of $\reals^n$ such that, for any $\xhat$ in 
$\reals^n\setminus S$, there exists a unit vector 
$v_{\xhat}$
such that, for any $x$ in $S$ and any minimal geodesic $\gamma ^*$
between $x= \gamma ^*(0)$ and $\xhat=\gamma ^*(\hat s)$,
with $\hat s>0$,
we have
$$
\frac{d\gamma ^*}{ds}(\hat s)^\top P(\xhat) \,  v_{\xhat}
\; <\; 0
\ .
$$
Then, $S$
is geodesically convex.
}
\end{lemma}
\begin{IEEEproof}
Assume that $S$ is not geodesically convex. Then, 
there is
a pair $(x_1,x_2) \in S$ such that, for any
minimal geodesic $\gamma _1^*$
between $x_1=\gamma _1^*(0)$ and $x_2=\gamma _1^*(s_2)$,
there exists $\hat s_1$ in $(0,s_2)$ for which
$\gamma _1^*(\hat s_1)$ is not in $S$.
Let
$
\xhat \;=\; \gamma _1^*(\hat s_1) \not \in S
$.
Note that
$
\gamma _2^*(s)\;=\; \gamma _1^*(s_2-s)
$
defines a minimal geodesic
between $x_2 = \gamma _2^*(0) \in S$ and
$\xhat =\gamma _2^*(\hat s_2)\not\in S$,
with $\hat s_2=s_2-\hat s_1 >0$.
With our assumption, 
since $x_1$ and $x_2$ are in $S$, there exists a unit vector $v_{\xhat}$
satisfying
$$
\frac{d\gamma _1^*}{ds}(\hat s_1)^\top P(\xhat) \,  v_{\xhat}
 \; <\; 0
\quad ,\qquad 
\frac{d\gamma _2^*}{ds}(\hat s_2)^\top P(\xhat) \,  v_{\xhat}
\; <\; 0
\ .
$$
But this impossible since we have
$
\frac{d\gamma _1^*}{ds}(\hat s_1)\;=\; -\frac{d\gamma _2^*}{ds}(\hat s_2).
$
\end{IEEEproof}

For Example~\ref{ex:1}, we shall see in the following section that, with the help of item 2a of 
Proposition~\ref{prop2}, 
for any $y$, the level set $\mathfrak{H}(y) = \{(x_1, x_2)\ :\ x_1 = y\}$ is geodesically convex for
the Riemannian metric given in \eqref{eqn:PexampleInx}.

As announced above, we conclude from Lemma \ref{lem2} that 
geodesic convexity of the levels sets of the output function is a 
necessary property in the ``general situation'' where (\ref{LP44})
holds (and when (\ref{LP7}) is an equality).
Actually, it is necessary, without any extra condition,
when the observer has an infinite gain margin.

\begin{definition}[infinite gain margin]
{\it
The observer $\dot\xhat = F(\xhat,y)$
for $\dot{x} = f(x)$
 is said to have an infinite 
gain margin with respect to $P$ if
(\ref{eqn:FandfEquality}) holds for every $x \in \reals^n$ and, for any geodesic $\gamma ^*$ 
minimal on $[0,\hat s)$, we have
\begin{equation}
\label{LP73}
\frac{d \gamma ^*}{ds}(s) P(\gamma ^*(s))
\left[F(\gamma ^*(s),h(\gamma ^*(0))-f(\gamma ^*(s))\right]\; <\; 0
\end{equation}
}
for all $s\in (0,\hat s)$.
\end{definition}

The term \textit{infinite gain margin} follows from the fact 
that,
if the observer $\dot\xhat = F(\xhat,y)$ makes $t\mapsto d(\hat 
X((\hat x,\x),t),X(x,t))$ 
nonincreasing (for each solution) and (\ref{LP73}) holds, then the same holds 
for the observer $\dot \xhat = f(\xhat)\;+\; \ell\,  
\left[F(\xhat,y)-f(\xhat)\right]$ for any real number $\ell > 1$.

\subsection{
Necessity of 
Uniform Detectability}
\label{sec1d}

The necessary condition in 
\eqref{LP71}
is linked to an observability property of the family of linear time-varying systems
obtained from linearizing \eqref{eqn:Plant1} along its solutions.
Assuming
the system
\eqref{eqn:Plant1} is 
forward
complete, for each $x$,
the corresponding solution to \eqref{eqn:Plant1} 
$t\mapsto \XX(x,t)$ is defined
on $[0,+\infty )$.
For each $x$,
the linearization of $f$ and $h$ evaluated along a solution 
$\XX(x,t)$
gives the
following functions defined on $[0,+\infty )$
$$
A_{x}(t)  \;=\;   \frac{\partial f}{\partial \x}(\XX(x,t))
\  ,\quad 
C_{x}(t) \;=\; \frac{\partial h}{\partial \x}(\XX(x,t)).
$$
These functions
define the following family of linear time-varying systems
with state $\xi \in \reals^n$ and output $\eta \in \reals^m$:
\begin{equation}
\label{LP6}
\dot \xi \;=\;  A_{x}(t)\,  \xi
\  ,\quad 
\eta \;=\; C_{x}(t)\,  \xi.
\end{equation}
Systems \eqref{LP6} are parameterized by the initial condition $\x$ of the
chosen solution $\XX(x,t)$.

The following theorem establishes a
relationship between
a detectability property
of \eqref{LP6}
and the existence of a 
bounded away from zero, upper bounded 
symmetric
covariant two-tensor
whose Lie derivative satisfies 
\eqref{LP71}.

\begin{theorem}
\label{thm:NecessitySatisfaction}
{\it
Assume system \eqref{eqn:Plant1} is
forward complete and
that there exist a
$C^1$ 
symmetric
covariant two-tensor
$P:\reals^n\to\reals^{n\times n}$
and strictly positive real numbers $\plower$ and $\pupper$
satisfying \eqref{LP71} and
\begin{equation}
\label{eqn:Pbounds}
0 \; < \;  \plower\,  I \; \leq \;  P(\x)  \leq \;  \pupper \,  I,
\quad \forall x\in \reals^n.
\end{equation}
Then,
for each $x \in \reals^n$, 
there exists a continuous\footnote{
We do not ask the function $K_x$ to be bounded.
}
function
$t\in [0,+\infty )\to K_x(t)$
such that the origin of the linear time-varying system
\begin{equation}
\label{eqn:LinearizationClosedLoop}
\dot \xi \;=\; \left(A_{x}(t) - K_{x}(t) C_{x}(t)\right)\,  \xi
\end{equation}
is uniformly exponentially stable.
}
\end{theorem}
\begin{IEEEproof}
To any $x \in \reals^n$, we associate the functions
$\Pi_x:[0,+\infty) \to \reals^{n\times n}$,
$K_{x}:[0,+\infty) \to \reals^{n}$, and
$\mathcal{V}_{x}:\reals^n\times [0,+\infty) \to \reals$
defined as
\begin{equation}\label{eqn:DefnsTraj}
\begin{array}{ll}\displaystyle
\Pi _{x}(t)\;=\;  P(\XX(x,t)),\ \
\mathcal{V}_{x}(\xi ,t)\;=\;  \xi^\top \Pi _{x}(t) \xi
\ ,\\ \displaystyle
K_{x}(t) \;=\;  \frac{\rho (\XX(x,t))}{2} \,  \Pi _{x}(t)^{-1}\,  
C_{x}(t)^\top
\  .
\end{array}
\end{equation}
We have
\begin{equation}
\label{LP24}
\plower\,  |\xi|^2 \; \leq \; \mathcal{V}_{x}(\xi ,t)\; \leq \; 
\pupper\,  |\xi|^2
\qquad \forall (x,t,\xi)
\end{equation}
and, with \eqref{LP71}, \eqref{eqn:D2dcurveBound}, \eqref{LP23}, and the definitions in \eqref{eqn:DefnsTraj}, we get
\begin{eqnarray*}
\frac{d}{dt} \left(v^\top \Pi _{x}(t)\,  v\right)
&= &
\left.\frac{\partial }{\partial 
\chi
}\left(v^\top
P(
\chi
)\,  v\right)\,  f(
\chi
)\right|_{
\chi
=\XX(x,t)}
\  ,
\\
&\leq&
-
 \frac{\varepsilon}{2}
v^\top \Pi _{x}(t)\,  v\\
& &
\;-\; 2\,  v^\top \Pi _{x}(t) \left(A_{x}(t) -K_{x}(t) C_{x}(t) \right)v
\  .
\end{eqnarray*}
Then, with (\ref{eqn:LinearizationClosedLoop}), 
we have
$
\frac{d}{dt} \mathcal{V}_{x}(\xi,t)
\; \leq \; -
 \frac{\varepsilon}{2}
\mathcal{V}_{x}(\xi,t)
$.
The conclusion follows with (\ref{LP24}).
\end{IEEEproof}

It follows from this proof that, if we do not have the upper bound 
$\pupper$ in (\ref{eqn:Pbounds}), we still have exponential stability, but we 
loose the uniformity property.
This would be the case, for instance, 
for the system (\ref{LP64})
of Example~\ref{ex:1} with 
$P$ given by \eqref{eqn:PexampleInx}
whose eigenvalues satisfy
\begin{eqnarray}
\label{LP67}
\begin{array}{ccl}
\lambda _{\min}(P(x))& \geq &\displaystyle \frac{(2+x_2^2)(1+x_1^2)- (x_1x_2-1 )^2}{3+x_2^2+x_1^2} = \\ 
& &\qquad \displaystyle
\frac{1+x_1^2+(x_1+x_2)^2}{3+x_2^2+x_1^2}\; \geq \; \frac{1}{3},
\end{array}
\end{eqnarray}
\begin{eqnarray}
\begin{array}{ccl}
\lambda _{\max}(P(x)) & \leq & \displaystyle
3+x_2^2+x_1^2.
\end{array}
\end{eqnarray}

\par\vspace{1em}
Exponential stability of the origin of (\ref{eqn:LinearizationClosedLoop})
is a detectability property for
\eqref{LP6}.
The necessity of this property for the existence of $P$
can be exploited to actually construct it,
as it will be shown in the companion paper.

\section{A Sufficient Condition}
\label{sec3}

In the previous section,
we assumed the existence of an observer 
making the function 
$t\mapsto d(\hat X((\hat x,\x),t),X(x,t))$ 
nonincreasing (respectively, strictly decreasing)
with $d$  being 
the distance associated with a Riemannian metric $P$.
We showed that $P$ has to satisfy
a (respectively, strict) inequality involving 
the output function. In this section,
we start from the data of such a
metric
 and investigate
the possibility of designing an observer 
making the
corresponding
Riemannian distance $d(\xhat,x)$
strictly decreasing along solutions.

In view of Theorem~\ref{thm3}, we assume that $P$ satisfies
$$
\mathcal{L}_f P(x)
\; \leq  \; 
\rho (x)\,  
\frac{\partial \h}{\partial \x}(\x) ^\top   \frac{\partial \h}{\partial \x}(\x) 
-\qlower\,  P(x)
\qquad \forall x\in\reals^n
$$
with $\qlower$ a strictly positive real number.
But, also, willing to be in a ``general situation'' in which (\ref{LP44}) holds 
and motivated by Lemma~\ref{lem2}, we restrict our attention
to the case where
the level set of the 
output function 
$\mathfrak{H}(y)$
is
geodesically convex for any $y$ in 
$\reals^m$.
Actually, we ask for the stronger (see Proposition \ref{prop2})
property that the sets $\mathfrak{H}(y)$ are totally geodesic
(see \cite[Section V.II]{Cartan.51}).

\begin{definition}[totally geodesic set]
{\it
Given a $C^1$ function $\varphi:\reals^n\mapsto \reals^m$
and a closed subset $\mathcal{C}$ of $\reals^n$,
the set
$$
S\;=\; \{x\in
\reals^n
:\,  \varphi(x)=0\}
\cap \mathcal{C}
$$
is said to be totally geodesic if, for any pair $(x,v)$ in 
$S\times\reals ^n$ such that
$$
\frac{\partial \varphi}{\partial x}(x)\,  v\;=\; 0
\quad ,\qquad
v^\top P(x)\,  v\;=\; 1
\  ,
$$
any geodesic $\gamma $ with
$$
\gamma (0)\;=\; x\quad ,\qquad 
\frac{d\gamma }{ds}(0)\;=\; v
$$
satisfies
$$
\varphi(\gamma (s)) \;=\; 0 \qquad \forall s \in J_\gamma 
\  ,
$$
where $J_\gamma $ is the maximal interval containing $0$ so that $\gamma (J_\gamma )$ 
is contained in $\mathcal{C}$.
}
\end{definition}
In the appendix, we establish a necessary and sufficient checkable condition
for the sets $\mathfrak{H}(y)$ to be totally geodesic.

\begin{example}[Motivational example -- continued]
For the system in Example~\ref{ex:1}, it is 
sufficient to check that the Christoffel symbol $\Gamma _{22}^1$ (see 
(\ref{eqn:ChristoffelSymbols}))
associated with 
the particular choice of $P$ in \eqref{eqn:PexampleInx}
for the family \eqref{eqn:PfamilyExample} 
is zero. In fact, we have
$
\Gamma _{22}^1 = \frac{1}{1+x_1^2+(x_1+x_2)^2}
\left(\begin{array}{@{\,  }c@{\quad }c@{\,  }}
1+x_1^2 & 1-x_1x_2 
\end{array}\right)
\left(\renewcommand{\arraystretch}{0.6}
\begin{array}{@{\,  }c@{\,  }}
0
\\
0
\end{array}\right) = 0
$.
\hfill$\Box$

\end{example}
\par\vspace{1em}
The following theorem gives a sufficient condition for the 
existence of an observer for the single output case.

\begin{theorem}
\label{thm1}
{\it
Assume there exist a complete {\color{black}$C^2$} Riemannian metric $P$ and
a set $\mathcal{C} \subset \reals^n$ such that
\begin{enumerate}
\item[H1~:]
$\mathcal{C}$ 
is geodesically convex,
closed, and 
with nonempty interior;
\item[H2~:]
there exist
a {\color{black}$C^1$} function $\rho :\reals^n\to \realsgeq$
and a strictly positive real number $q$ such that
\begin{equation}
\label{LP10}
\mathcal{L}_f P(x)
\; \leq  \; 
\rho (x)\,  
\frac{\partial \h}{\partial \x}(\x) ^\top   \frac{\partial \h}{\partial \x}(\x) 
-\qlower\,  P(x)
\qquad \forall x\in\mathcal{C},
\end{equation}
\item[H3~:]
The number of outputs is $m=1$ and,
for each $y$ in $h(\mathcal{C})$,
the set $\mathfrak{H}(y)
\cap \mathcal{C}$
is totally geodesic.
\end{enumerate}
Then,
 for any positive real number $E$ there exists 
a continuous function $k_E:\reals^n\to \reals$ such that,
with the observer given by
\begin{equation}
\label{LP9}
F(\xhat,y)\;=\; f(\xhat) \;-\; k_E(\xhat)\,  P(\xhat)^{-1}
\frac{\partial h}{\partial x}(\xhat) ^\top
\frac{\partial \delta }{\partial y_1}(\h(\xhat),y)
\  ,
\end{equation}
where
\begin{equation}
\label{LP68}
\delta (y_1,y_2)\;=\; |y_1-y_2|^2
\  ,
\end{equation}
the following holds
(see (\ref{LP57})):
\begin{equation}
\label{LP40}
\begin{array}{lll}
\Did d(\xhat,\x)
& \leq & \displaystyle -
\frac{\qlower}{4} 
\,  d(\xhat,\x)  \\
& & \displaystyle
\hspace{-0.8in} \forall (x,\xhat) \in
\left\{(x,\xhat )\,  :\:
d(\xhat,x) <  E
\right\}
\; \bigcap\; 
\left(\mbox{\rm int}(\mathcal{C}) \times \mbox{\rm int}(\mathcal{C})\right)
\  .
\end{array}
\end{equation}
Moreover, expression (\ref{LP9}) is intrinsic (i.e., coordinate 
independent) and gives an observer with infinite gain margin.
}
\end{theorem}

\begin{example}[Motivational example -- continued]
\label{ex:1-rere}
We have 
already checked that, for the system
(\ref{LP64}) and with $P$ given in \eqref{eqn:PexampleInx}
 all the conditions of Theorem~\ref{thm1} hold 
globally, i.e.,~with $\mathcal{C}=\reals^2$. 
 Hence, the observer
given by \eqref{LP9} becomes
\begin{eqnarray*}
\left(\begin{array}{lll}
\dot{\hat x}_1 \\ \dot{\hat x}_2
\end{array}\right) & = & \left(\begin{array}{c}
\hat x_2 \sqrt{1+\hat x_1^2}\\ -\displaystyle\frac{\hat x_1\hat x_2^2}{\sqrt{1+\hat x_1^2}} 
\end{array}\right) \\
& & \hspace{-0.3in}\;-\; \frac{2k_E(\xhat) }{1+\xhat _1^2+(\xhat _1+\xhat_2)^2}
 \left(\begin{array}{c}
1+\xhat _1^2 \\ 1-\xhat _1\xhat _2
\end{array}\right)
(\xhat _1 -y)
\  .
\end{eqnarray*}
\null \hfill $\Box$\\[1em]
\end{example}

\begin{remark}~
\normalfont
\begin{itemize}
\item
Theorem~\ref{thm1}
gives
a (nonglobal) solution to problem \OurProblem.
When the assumptions of Theorem~\ref{thm1} hold 
globally, i.e., they hold for $\mathcal{C}=\reals^n$, 
the observer given by 
(\ref{LP9}) guarantees convergence of 
the estimated state to the system state, 
semiglobally with respect to the zero estimation error set $\A$.

The fact that we do not get global asymptotic stability 
is likely 
due to the elementary form of the observer (\ref{LP9})
and its infinite gain margin.
We expect that other choices for this observer are possible to obtain a global asymptotic stability result.
\item
As discussed in \ref{sec1b}, 
we do not claim
in Theorem~\ref{thm1}
that the flow generated 
by the observer has a contraction property but simply that the 
Riemannian distance between estimated state and system state decays 
along the solutions.
In other words, this result establishes that 
the function $(\xhat,x)\mapsto d(\xhat,x)$ can be used as a Lyapunov function
for the zero error set $\mathcal{A}$ and guarantees this function has 
an exponential decay along the solutions. But it does no say that
$d(\xhat_1,\xhat_2)$ decays along two arbitrary solutions of 
the flow generated by the observer.\hfill$\Box$
\end{itemize}
\end{remark}
\par\vspace{1em}

Theorem~\ref{thm1} is a direct consequence of the following lemma (for which there
is no restriction on the number of outputs)
and the fact that, when the 
number of outputs is $m=1$, assumption H3 implies
the assumption H3' of the lemma; see Proposition~\ref{prop2}.

\begin{lemma}
\label{lem3}
{\it 
Assume there exist a complete {\color{black}$C^2$} Riemannian metric
$P$,
a set $\mathcal{C} \subset \reals^n$,
a {\color{black}$C^1$} function $\rho :\reals^n\to \realsgeq$,
and a strictly positive real
number $q$ satisfying H1 and H2
 of Theorem~\ref{thm1}.
Assume also
there exists
a {\color{black}$C^3$} function $\delta :\reals^m\times\reals^m\to\realsgeq$
satisfying
\begin{equation}
\label{LP58}
\delta (h(x),h(x))\;=\; 0\quad ,\qquad 
\left.
\frac{\partial ^2\delta }{\partial y_1^2}(y_1,y_2)
\right|_{y_1=y_2=h(x)}\; >\; 0
\end{equation}
for all $x\in 
\mathcal{C}$,
and, such that
\begin{enumerate}
\item[H3':]
for any pair $(x_1,x_2)$ in $\mathcal{C}\times 
\mathcal{C}$ satisfying
$$
h(x_1)\;\neq \; h(x_2)
$$
and for any minimal geodesic $\gamma ^*$ between
$x_1=\gamma ^*(s_1)$ and $x_2=\gamma ^*(s_2)$
satisfying
$\gamma ^*(s) \in \mathcal{C}$ for all $s \in [s_1, s_2]$,
with $s_1\leq s_2$,
we have
\begin{equation}
\label{LP39}
\frac{d}{ds}\delta (h(\gamma ^*(s)),h(\gamma ^*(s_1)))\; >\; 0
\qquad \forall s\in (s_1,s_2]
\  .
\end{equation}
\end{enumerate}
Then,
 the claim of Theorem~\ref{thm1} holds true with a function 
 $\delta $ satisfying H3' (instead of
 $\delta$ as in (\ref{LP68})).
 }
\end{lemma}

\begin{remark}~
\normalfont
\begin{itemize}
\item
Property H3' says that we can find a ``distance-like'' function 
$\delta $ in the output space allowing us to express that the output function $h$ 
preserves some kind of monotonicity. Namely, as  the 
distance increases along a geodesic in  the state space, the same holds in the output 
space measured by $\delta $.
This property has some relationship with the notions of metric-monotone function
introduced in \cite{Price.40} and of
geodesically monotone function defined in \cite[Definition 6.2.3]{Rapcsak.97}.
In the appendix, we establish a connection with totally 
geodesic sets and geodesic convexity. 

With such a property, by following a descent 
direction for the ``distance'' in the output space, we are guaranteed 
to decrease the distance in the state space. This feature
is
exploited in the observer given by (\ref{LP9}) via a  high-gain 
term which enforces that such a descent direction is dominating.
\item
Property H3' with $\delta (y_1,y_2)=|y_1-y_2|^2$
has been invoked already in \cite{Tsinias.90} but for 
the case when $P$ is constant.\hfill$\Box$
\end{itemize}
\end{remark}

\begin{IEEEproof}
Note that since we have
$$
\xhat\;=\; x
\qquad \Rightarrow\qquad 
F(\hat x,y)\;=\; f(\xhat)\;=\; f(x),
$$
the result already holds when $d(x,\xhat)$ is zero. 
Therefore, the remainder of the proof 
only considers pairs $(\xhat ,x)$ that are in
\mbox{$\left(\mathcal{C}\times\mathcal{C}\right)\setminus\A$}.

The Riemannian metric $P$ being complete,
any geodesic is defined on $(-\infty ,+\infty )$ and
the Riemannian distance
$d(x_1,x_2)$ is given by the 
length of a minimal geodesic $\gamma ^*$ between $x_1$ and $x_2$.
Since $\mathcal{C}$ is geodesically convex by H1,
for any pair $(x_1,x_2)$ in
\mbox{$\left(\mathcal{C}\times\mathcal{C}\right)\setminus\A$},
there exists a minimal geodesic $\gamma ^*$
between $x_1 =\gamma ^*(s_1)$ and $x_2=\gamma ^*(s_2)$ satisfying
$
\gamma ^*(s) \in \mathcal{C}$ for all $s \in [s_1,s_2]$.

Let $(\xhat ,x)$ be any pair  in
\mbox{$\left(\mathcal{C}\times\mathcal{C}\right)\setminus\A$}
and
$\gamma ^*$ denote a minimal geodesic
between $x=\gamma ^*(0)$ and $\xhat=\gamma ^*(\hat s)$
satisfying
$
\gamma ^*(s) \in \mathcal{C}$ for all $s \in [0,\hat s]$.  With $y=h(x)$, 
take $F$ as in (48). 
It gives
$$
\displaylines{\quad 
\frac{d\gamma ^*}{ds}(\hat s)^\top P(\gamma ^*(\hat s))  
\left[\F(\gamma ^*(\hat s),\y) - f(\gamma ^*(\hat s))\right]\hfill
\cr\hfill
- \frac{d\gamma ^*}{ds}(0)^\top P(\gamma ^*(0))
\left[\F(\gamma ^*(0),\y)-f(\gamma ^*(0))\right]
\hfill
\cr\hfill
\vbox{\ialign{\strut$\displaystyle{#}$\hfill&$\quad$ 
\hfill$\displaystyle{#}$\crcr
=\;
-k_E(\xhat)\, 
\frac{d\,  h\circ \gamma ^*}{ds}(\hat s)^{\top}
\frac{\partial \delta }{\partial y_1}(h(\gamma ^*(\hat s)),y)
.
&\refstepcounter{equation}\label{LP35}
(\theequation)\cr
\crcr }}\cr}$$
On the other hand, we have
\\[0.7em]$\displaystyle 
\frac{d\gamma ^*}{ds}(\hat s)^\top P(\xhat)  \,  f(\xhat)
- \frac{d\gamma ^*}{ds}(0)^\top P(x)\,  f(x)
$\refstepcounter{equation}\label{Aux1}\hfill$(\theequation)$
\\[0.5em]\null\hfill$\displaystyle 
=\; 
\int_0^{\hat s}\frac{d}{ds} \left(\frac{d\gamma ^*}{ds}(s)^\top P(\gamma ^*(s))  
\,  f(\gamma ^*(s))\right) ds
\  .
$\\[0.7em]
Also the Euler-Lagrange form of the geodesic equation
reads, for the
$i$-th coordinate,
\begin{equation}\non
\begin{array}{lll}\displaystyle
2\,  \frac{d}{ds} \left(\sum_kP_{ik}(\gamma ^*(s)) \frac{d\gamma _k
^*}{ds}(s)\right)
& = & \\
\displaystyle
& &\displaystyle \hspace{-0.9in} \sum_{k,l}\frac{d\gamma _k^*}{ds}(s)^\top\,  
\frac{\partial P_{kl}}{\partial x_i}(\gamma ^*(s))
\,  {\frac{d\gamma _l^*}{ds}(s)}
\  .
\end{array}
\end{equation}
Then, with the definition of the Lie 
derivative $\mathcal{L}_fP$ and (47), 
we get
$$\displaylines{\quad
\frac{d}{ds} \left(\frac{d\gamma ^*}{ds}(s)^\top P(\gamma ^*(s))  
\,  f(\gamma ^*(s))\right)
\hfill\null\cr\hfill
\vbox{\ialign{\strut$\displaystyle{#}$\quad  \null \crcr
\qquad \qquad \qquad \qquad =\; 
\frac{1}{2}\,  
\frac{d\gamma ^*}{ds}(s)^\top \mathcal{L}_fP(\gamma ^*(s))  
\frac{d\gamma ^*}{ds}(s)
\  ,
\cr
\qquad \qquad \qquad \qquad\leq \; 
\frac{\rho (\gamma ^*(s))}{2}
\left|
\frac{\partial \h}{\partial \x}(\gamma ^*(s)) 
\frac{d \gamma ^*}{ds}(s)\right|^2 \cr
\qquad \qquad \qquad \qquad \qquad\;-\; 
\frac{\qlower}{2}
\, 
\frac{d\gamma ^*}{ds}(s)^\top  P(\gamma ^*(s)) 
\frac{d\gamma ^*}{ds}(s)
\cr
\qquad \qquad \qquad \qquad \leq\; 
\frac{\rho (\gamma ^*(s))}{2}
\left|\frac{d\,  h \circ \gamma ^*}{ds}(s)\right|^2
\;-\;
\frac{\qlower}{2}
\  ,
\cr
\crcr 
}}
\refstepcounter{equation}\label{Aux2}
\hfill(\theequation)
}
$$
where, in the last inequality, we have used
$$
\frac{d\gamma ^*}{ds}(s)^\top  P(\gamma ^*(s)) 
\frac{d\gamma ^*}{ds}(s)=1
$$
since $\gamma^*$ is normalized.
With
$
d(\xhat,x)\;=\; 
\hat s
$
as given in (12), 
replacing \eqref{Aux2} into \eqref{Aux1}
yields
\\[0.7em]$\displaystyle 
\frac{d\gamma ^*}{ds}(\hat s)^\top P(\gamma ^*(\hat s))  \,  f(\gamma ^*(\hat s))
- \frac{d\gamma ^*}{ds}(0)^\top P(\gamma ^*(0))\,  f(\gamma ^*(0))
$\\[0.7em]\null\hfill$\displaystyle 
\qquad\;\leq \; 
\int_0^{\hat s}
\frac{\rho (\gamma ^*(s))}{2}
\left|\frac{d\,  h \circ \gamma ^*}{ds}(s)\right|^2 ds
\;-\;
\frac{\qlower}{2}\,  d(\xhat,x)
\  .$\refstepcounter{equation}\label{Aux01}
\hfill$(\theequation)$\\[0.7em]
Then, from (36), 
using (\ref{LP35}) and \eqref{Aux01}, we obtain
$$\displaylines{
\Did d(\xhat,\x)
\hfill\null\cr\hfill
\vbox{\ialign{\strut$\displaystyle{#}$\crcr
\begin{array}{@{}cl@{}}
\leq& \displaystyle
\left[\frac{d\gamma ^*}{ds}(\hat s)^\top P(\gamma ^*(\hat s))  
\left(\F(\gamma ^*(\hat s),\y) - f(\gamma ^*(\hat s))\right)
\right.
\\[0.5em]
&\displaystyle
\left.- \frac{d\gamma ^*}{ds}(0)^\top P(\gamma ^*(0))
\left(\F(\gamma ^*(0),\y)-f(\gamma ^*(0))\right)
\right]
\\[0.5em] 
&\displaystyle
\hspace{-0.3in}+ \left[\frac{d\gamma ^*}{ds}(\hat s)^\top P(\gamma ^*(\hat s))  
f(\gamma ^*(\hat s))
- \frac{d\gamma ^*}{ds}(0)^\top P(\gamma ^*(0))
f(\gamma ^*(0))\right]
\end{array}
\cr
\crcr }}\cr}$$
\\[0.5em]$\displaystyle 
\leq \; 
-\;  k_E(\xhat)\,  
\frac{d\,  h\circ \gamma ^*}{ds}(\hat s)^{\top}
\frac{\partial \delta }{\partial y_1}(h(\gamma ^*(\hat s)),y)^\top
$\refstepcounter{equation}\label{LP59}\hfill$(\theequation)$
\\[0.7em]\null\hfill$\displaystyle 
+ \int_0^{\hat s}
\frac{\rho (\gamma ^*(s))}{2}
\left|\frac{d\,  h \circ \gamma ^*}{ds}(s)\right|^2
ds
\: -
\frac{\qlower}{2}\,  d(\xhat,x)
\:  .
$

To proceed it is appropriate to associate two functions $a$ and $b$ 
to any triple $(\xhat,x,\gamma ^*)$ with
$(\xhat ,x)$ in
\mbox{$\left(\mathcal{C}\times\mathcal{C}\right)\setminus\A$}
and
$\gamma ^*$, a minimal geodesic
between $x=\gamma ^*(0)$ and $\xhat=\gamma ^*(\hat s)$
satisfying
$
\gamma ^*(s) \in \mathcal{C}$ for all $s \in [0,\hat s]$. These 
functions are defined on $[0,\hat s]$ as follows:\footnote{
When $\hat s =0$ the functions $a_{(\xhat,x,\gamma ^*)}$
and $b_{(\xhat,x,\gamma ^*)}$ are only defined at zero.
}
\begin{eqnarray*}
a_{(\xhat,x,\gamma ^*)}(r)&=&\frac{1}{r}\,  
\frac{d\,  h\circ \gamma ^*}{ds}(r)^{\top}
\frac{\partial \delta }{\partial y_1}(h(\gamma ^*(r)),h(\gamma ^*(0)))^\top
\end{eqnarray*}
if  $0 < r \leq \hat s$, and
\begin{eqnarray*}
a_{(\xhat,x,\gamma ^*)}(0) &=& \\
& & \hspace{-0.6in}
\frac{d\,  h\circ \gamma ^*}{ds}(0)^\top
\frac{\partial ^2 \delta }{\partial y_1^2}(h(\gamma ^*(0)),h(\gamma ^*(0)))^\top\frac{d\,  h\circ \gamma ^*}{ds}(0);
\end{eqnarray*}
and
\begin{eqnarray*}
b_{(\xhat,x,\gamma ^*)}(r)&=&
\frac{1}{r}\,  
\int_0^{r}
\frac{\rho (\gamma ^*(s))}{2}
\left|\frac{d\,  h \circ \gamma ^*}{ds}(s)\right|^2
ds
\end{eqnarray*}
if  $0 < r \leq \hat s$, and
\begin{eqnarray*}
b_{(\xhat,x,\gamma ^*)}(0)&=&
\frac{\rho (\gamma ^*(0))}{2}
\left|\frac{d\,  h \circ \gamma ^*}{ds}(0)\right|^2.
\end{eqnarray*}
We remark with (51) 
that $\delta $ reaches its global minimum 
at $y_1=y_2=h(x)$. This implies
\begin{eqnarray*}
\frac{\partial \delta }{\partial y_1}(h(\gamma ^*(r)),h(\gamma ^*(0)))
& = & \\
& & 
\hspace{-1.5in}
\left[\int_0^1
\left(\frac{\partial ^2\delta }{\partial y_1^2}(h(\gamma ^*(\sigma r)),\gamma ^*(0))
\frac{d\,  h \circ \gamma ^*}{ds}(\sigma r)\right)
d\sigma \right] r
\end{eqnarray*}
for all $r\in [0,\hat s]$.
As a consequence, the functions $a$ an $b$ are continuous on
$[0,\hat s]$. Moreover the property H3' gives readily the implication
$$
h(x)\;\neq \; h(\xhat)\qquad \Longrightarrow\qquad 
a_{(\xhat,x,\gamma ^*)}(r)\; >\; 0\quad \forall r\in (0,\hat  s]
\  .
$$
In the case when $h(x)\;= \; h(\xhat)$, we are only left with the following two 
possibilities:
\begin{enumerate}
\item
$h\circ \gamma ^*$ is constant on $[0,\hat s]$. Then we have
$
\frac{d h\circ \gamma ^*}{ds}(s)\;=\; 0$ for all $s\in[0,\hat s]$
and therefore
$
a_{(\xhat,x,\gamma ^*)}(r)\;=\; b_{(\xhat,x,\gamma ^*)}(r)\;=\; 0$ for all $r\in[0,\hat s].
$
\item
$h\circ \gamma ^*$ is not constant on $[0,\hat s]$. Then, there exists 
some $s_1$ in $(0, \hat s]$ such that
$
h(\gamma (s_1))\; \neq \; h(\gamma ^*(0))=h(x).
$
With H3', this implies that the function $s\mapsto \delta (h(\gamma ^*(s)),h(\gamma ^*(0)))$ 
is not constant on $[0, \hat s]$. But since we have
$
\delta (h(\gamma ^*(\hat s)),h(\gamma ^*(0)))= 
\delta (h(\gamma ^*(0)),h(\gamma ^*(0)))= 0,
$
this function must
reach a maximum at some point $s_m$ in $(0,\hat s)$ where we have
$$
\delta (h(\gamma ^*(s_m)),h(\gamma ^*(0))) >0,$$
$$\frac{d}{ds}\delta (h(\gamma ^*(s_m)),h(\gamma ^*(0)))= 0,
$$
and therefore
$
h(\gamma ^*(s_m))\; \neq \; h(\gamma ^*(0)).
$
But this contradicts H3'. So this case is impossible.
\end{enumerate}
In any case, we have established that $a_{(\xhat,x,\gamma ^*)}(\hat 
s)$ is non negative and if it is zero then
$
b_{(\xhat,x,\gamma ^*)}(r)= 0$ for all $r\in [0,\hat s].
$

Now,
let $\check x$ be an arbitrary point in $\mathcal{C}$. Call it origin.
For each integer $i$, we introduce the set
$$
\mathcal{K}_i
\;=\;
 \{(x,\xhat )
\in \mathcal{C}\times \mathcal{C}
\,  :\: d(\xhat,x)\leq   
E\; ,\  i\leq 
d(\check x ,\xhat)
\leq i+1\}
\  .
$$
From the  Hopf-Rinow Theorem 
\cite[Theorem II.1.1]{Sakai.96} 
$\mathcal{K}_i$ is compact.

To conclude it is sufficient to
prove the existence of a real number $k_i$ such 
that, for any pair $(\xhat ,x)$ in
\mbox{$\mathcal{K}_i\setminus\A$}
and any minimal geodesic $\gamma ^*$
between $x=\gamma ^*(0)$ and $\xhat =\gamma ^*(\hat s)$
satisfying
$
\gamma ^*(s)\in\mathcal{C}$ for all $s\in [0,\hat s],
$
we have
\\[0.5em]\null \hfill $\displaystyle 
\frac{q}{4}\;+\; k_i\,  a_{(\xhat,x,\gamma ^*)}(\hat s) > 
b_{(\xhat,x,\gamma ^*)}(\hat s).
$\hfill \null \\[0.5em]
Indeed, with this inequality, the definitions of $a$ and $b$ and
(\ref{LP59}) where $d(\xhat,x)=\hat s$,
we obtain (50) 
provided the function $k_E$ satisfies
\\[0.3em]\null \hfill 
{\color{black}$\displaystyle 
k_E(\xhat)\; \geq \; k_i
\qquad \forall \xhat
\in\mathcal{C}
:\, i\leq 
d(\xhat ,x)
\leq i+1
\  .
$\hfill \null
}
\par

Proceeding by contradiction, suppose that such $k_i$ does not exist.
Then, there exists a sequence $(\hat s_n, x_n,\hat 
x_n, \gamma _n^*)$, with 
$\hat s_n \geq 0$,
$(x_n,\hat x_n)$ in
\mbox{$\mathcal{K}_i\setminus\A$},
and $\gamma _n^*$ a minimal geodesic
between $x_n=\gamma _n^*(0)$ and $\xhat _n=\gamma _n^*(\hat s_n)$ 
satisfying
$
\gamma _n^*(s)\in\mathcal{C}$ for all $s\in [0,\hat s_n]$
and
\\[0.5em]\null \hfill $\displaystyle 
\frac{q}{4}\;+\; n\,  a_{(\xhat _n,x_n,\gamma _n^*)}(\hat s_n)\; \leq \; 
b_{(\xhat _n,x_n,\gamma _n^*)}(\hat s_n)
\  .
$\refstepcounter{equation}\label{LP16}\hfill(\theequation)
\\[0.7em]
{\color{black}
Moreover, the functions $a_{(\xhat,x,\gamma ^*)}$ and $b_{(\xhat,x,\gamma 
^*)}$ are $C^1$ on $[0,\hat s]$. Indeed, they can be written as
$$
a_{(\xhat,x,\gamma ^*)}=\frac{f_a(r)}{r}
\  ,\quad 
a_{(\xhat,x,\gamma ^*)}=\frac{f_b(r)}{r}
\forall r \in ]0,\hat s]
$$
where the function $f_a$, respectively $f_b$, is $C^2$ since
$h$, $\gamma ^*$ and $\delta $ are $C^3$,
respectively, $\rho $ 
is $C^1$ and $h$ and $\gamma ^*$ are $C^2$. 
We have the following technical property.
\par\vspace{1em}\noindent
\textit{%
Claim 1:
Let $f$ be a $C^2$ function defined on a neighborhood of $0$ in 
$\RR$, where it is $0$.
The function $\varphi $ defined as $\varphi (r)=\frac{f(r)}{r}$ 
if $r\neq 0$ and $\varphi (0)=f^\prime(0)$ is $C^1$.}
\par\vspace{1em}\noindent
\textit{Proof~:}
Clearly, $\varphi $ is $C^2$ everywhere except may be at $0$. Its first 
derivative is $\varphi ^\prime(r)=\frac{f(r)-rf^\prime(0)}{r^2}$. 
It is also continuous at $0$ since $\lim_{r\to 0}\varphi 
(r)=f^\prime(0)=\varphi (0)$. Its first derivative at $0$ exists if
$\lim_{r \to 0}\frac{\varphi (r)-\varphi (0)}{r}=\lim_{r\to 
0}\frac{f(r)-rf^\prime(0)}{r^2}$ exists, which is the case since, 
due to $f$ being $C^2$, we have
\begin{eqnarray*}
\frac{f(r)-rf^\prime(0)}{r^2}&=&\frac{1}{r^2}\int_0^r 
[f^\prime(s)-f^\prime(0)]ds
\\
&=&\frac{1}{r^2}\int_0^r\int_0^s 
f^{\prime\prime}(t) dt ds
\\
&=&
\frac{1}{r^2}\int_0^r f^{\prime\prime}(t) [r-t]dt
\end{eqnarray*}
which leads to
$\varphi ^\prime(0)=\frac{1}{2}f^{\prime\prime}(0)$. We have also 
$$
\frac{f(r)-rf^\prime(r)}{r^2}=-\frac{1}{r^2}\int_0^r 
s f^{\prime\prime}(s) ds
$$
This implies
$$
\lim_{r\to 0} \varphi ^\prime(r) = \varphi^\prime(0)
$$
and therefore $\varphi ^\prime $ is continuous.
}

We also have the following claim.

\noindent
{\it Claim 2: There exists a subsequence $(\hat s_{n_1},x_{n_1},\xhat _{n_1},\gamma _{n_1}^*)$ of $(\hat s_{n},x_{n},\xhat _{n},\gamma _{n}^*)$ such that
\begin{eqnarray}\label{eqn:LimitingSolPoints}
\lim_{n_1\to \infty } (\hat s_{n_1},x_{n_1},\xhat _{n_1})\!\!\! &=& \!\!\!
(\hat s_\omega , x_\omega ,\xhat _\omega ),\\
\label{eqn:LimitingSol}
\lim_{n_1\to \infty }\gamma _{n_1}^*(s)\!\!\! &=&\!\!\! \gamma _\omega (s)
\ \textrm{uniformly in} \  s \in [0, E]
\  , \qquad
\end{eqnarray}
where $\gamma_\omega:[0,\hat s_\omega]\to \mathcal{C}$ is  a minimal geodesic between 
$x_\omega$ and $\hat x_\omega$.}

\noindent
To prove the claim, not that since $(x_n,\hat x_n)$ is in the compact set
$\mathcal{K}_i$ and $\gamma _n^*$ is a minimal geodesic
taking values in $\mathcal{C}$ when restricted to $[0,\hat s_n]$,
from
\\[0.7em]$\displaystyle 
\sqrt{\plower}\,  |x_1-x_2|\; \leq  \; d(x_1,x_2)\; \leq  \; 
\sqrt{\pupper}\,  |x_1-x_2|
$\hfill \null \\\null \hfill $\displaystyle 
\forall (x_1,x_2)\in \mathcal{C}\times\mathcal{C}
\  ,
$\\[0.7em]
we get
$$
\sqrt{\plower}
\,  |\gamma _n^*(s)-x_n|\: \leq \: d(\gamma _n^*(s),x_n)\: \leq \: \hat s_n
\: \leq  \: E
\  \quad \forall s\in [0,\hat s_n]
$$
and
\begin{eqnarray*}
|x_n|
\; \leq \;
|\xhat _n-x_n|\;+\; |\xhat _n|
&\leq &
\frac{d(\xhat _n,x_n)}{\sqrt{\plower}}\;+\; (i+1)
\  ,
\\
&\leq &
\frac{E}{\sqrt{\plower}}\;+\; (i+1)
\  .
\end{eqnarray*}
This implies that {\color{black}$\gamma _n^*:[0,E]\to \mathcal{C}$}
takes its values in a compact set
independent of the index $n$. Moreover, $\gamma _n^*$ being a solution of the geodesic 
equation,
there exists a subsequence with index 
$n_1$ and a quadruple $(\hat s_\omega , x_\omega ,\hat x_\omega , \gamma _\omega )$ such that \eqref{eqn:LimitingSolPoints}-\eqref{eqn:LimitingSol} hold (see, for instance,  \cite[Theorem 5, \S 1]{Filippov.88}),
where $\gamma _\omega $ is a solution of the geodesic equation and,
since $\mathcal{C}$ is closed, it satisfies
$$
\gamma _\omega (0) \;=\;  x_\omega
\  ,\quad 
\gamma _\omega (\hat s_\omega ) \;=\;  \xhat _\omega
\  ,\quad 
\gamma _\omega (s)\in \mathcal{C}\quad  \forall s\in [0,\hat s_\omega ]
\  .
$$
Finally, according to \cite[Lemma II|.4.2]{Sakai.96}, it is minimizing 
between $x_\omega$ and $\xhat _\omega$.

Now, the functions $h$,
$\rho $
and $\frac{\partial  h}{\partial x}$ restricted to 
the compact set where the functions $\gamma _n^*$ take their values,
are continuous and bounded.  Also, from the geodesic equation and completeness, 
the same holds for $\gamma _n^*$, $\frac{d\gamma _n^*}{ds}$ and 
$\frac{d^2\gamma _n^*}{ds^2}$ restricted to
$[0, \hat s_n]$. 
With the definition of $b_{(\xhat _n,x_n,\gamma _n^*)}$, this 
implies that the right-hand side of (\ref{LP16}) is upper bounded, say by 
$B$. Consequently, we have
$$
\frac{q}{4}\;+\; n\,  a_{(\xhat _n,x_n,\gamma _n^*)}(\hat s_n)\; \leq 
\; B \qquad \forall n.
$$
Since $a_{(\xhat _n,x_n,\gamma _n^*)}(\hat s_n)$ is nonnegative, 
this implies that
$
a_{(\xhat _\omega ,x_\omega ,\gamma _\omega )}(\hat s_\omega )
= 0.
$

{\color{black}
If $\xhat _\omega \neq x_\omega$, since
$a_{(\xhat _\omega ,x_\omega ,\gamma _\omega )}(\hat s_\omega )$ 
is zero, we have seen that the same holds for
$b_{(\xhat _\omega ,x_\omega ,\gamma _\omega )}(r)$, for all
for all $r\in [0,\hat s_\omega ]$.
}
On the other hand, (\ref{LP16}) yields
$$
\frac{q}{4}\; \leq \; 
b_{(\xhat _\omega ,x_\omega ,\gamma _\omega ^*)}(\hat s_\omega )
$$
where $q$ is strictly positive. So we have a contradiction.

If $\xhat _\omega = x_\omega$,
also by compactness, there exists a subsequence with index 
$n_2$ of the subsequence with index $n_1$ in Claim 1 such that we have
$$
v_\omega \;=\; \lim_{n_2 \to \infty }
\frac{\xhat _{n_2}-x_{n_2}}{d(\xhat _{n_2},x_{n_2})}
\;=\; \lim_{n_2 \to \infty }
\frac{\xhat _{n_2}-x_{n_2}}{\hat s _{n_2}}
$$
Note that since $\xhat _\omega = x_\omega$, we have $\hat s_{n_1}$ 
(and also $\hat s_{n_2}$) converging to zero.
But, with the identity
$$
\xhat _{n_2}\;=\; x_{n_2}\;+\; \int_0^{\hat s_{n_2}}
\frac{d\gamma^* _{n_2}}{ds}(s) ds
\  ,
$$
this gives also
$$
v_\omega \;=\; \lim_{n_2 \to \infty }
{\color{black}
\frac{\gamma^*_{n_2}(\hat s_{n_2}) - \gamma^*_{n_2}(0)}{\hat s_{n_2}}}
\;=\; \frac{d\gamma _\omega ^*}{ds}(0)
$$
{\color{black}
On the other hand, since the functions
$a_{(\xhat _n,x_n,\gamma _n^*)}$ and $b_{(\xhat _n,x_n,\gamma _n^*)}$
are $C^1$ on $[0,E]$, and the way they depend on $n$ is only via
$\gamma _n^*$ (which takes its values in a compact 
set independent of $n$), there exist real numbers $A_1$ and $B_1$ 
such that we have
\begin{eqnarray*}
a_{(\xhat _n,x_n,\gamma _n^*)}(\hat s_n)&\geq &
a_{(\xhat _n,x_n,\gamma _n^*)}(0)
\;-\; A_1\,  \hat s_n\  ,
\\
b_{(\xhat _n,x_n,\gamma _n^*)}(\hat s_n)&\leq &
b_{(\xhat _n,x_n,\gamma _n^*)}(0)
\;+\; B_1\,  \hat s_n\  .
\end{eqnarray*}
}
Since we are in the case where $\hat s_{n_1}$ goes to $0$, this implies
\begin{eqnarray*}
0&=&\lim_{n_2 \to \infty } a_{(\xhat _{n_2},x_{n_2},\gamma _{n_2}^*)}(0)
\  ,
\\
&=& v_\omega ^\top\,  \frac{\partial h}{\partial x}(x_\omega )^\top
\frac{\partial ^2 \delta }{\partial y_1^2}(x_\omega ,x_\omega )^\top
\frac{\partial h}{\partial x}(x_\omega )v_\omega 
\  .
\end{eqnarray*}
With (51), we obtain
$$
\frac{\partial h}{\partial x}(x_\omega )v_\omega \;=\; 0
$$
and therefore~:
\\[0.7em]$\displaystyle 
\lim_{n_2\to \infty }
b_{(\xhat _{n_2},x_{n_2},\gamma _{n_2}^*)}(\hat s_{n_2})
$\hfill \null \\\null \hfill $
\begin{array}{@{}rcl@{}}
\leq \; 
\displaystyle \lim_{n_2\to \infty }
b_{(\xhat _{n_2},x_{n_2},\gamma _{n_2}^*)}(0)&=&\displaystyle 
\frac{\rho (\gamma _\omega ^*(0))}{2}
\left|\frac{\partial h}{\partial x}(x_\omega )v_\omega \right|^2
\  ,
\\
&=& 0
\  .
\end{array}
$\\[0.7em]
This contradicts (\ref{LP16}).

So we have stablished the existence of $k_i$.

Finally, in (\ref{LP35}), we have, with (52), 
\begin{eqnarray*}
\frac{d\,  h\circ \gamma ^*}{ds}(\hat s)^{\top}
\frac{\partial \delta }{\partial y_1}(h(\gamma ^*(\hat s)),y) & & \\ 
& & \hspace{-0.8in} =
\frac{d}{d\hat s}\delta (h(\gamma ^*(\hat s)),h(\gamma ^*(s)))\; >\; 0
\end{eqnarray*}
and
$
\F(\gamma ^*(0),\y)\;=\; f(\gamma ^*(0))
$.
So (39) 
holds and the observer has an infinite gain margin.

To prove the last point of Theorem 3.3, 
let
$\phi$ define a diffeomorphism as in (33). 
Let $\bar h$, $\bar k_E$, $\bar f$, $\bar F$ and $\bar P$  
be the expressions of 
$h$, $k_E$, $f$, $F$ and $P$
respectively in the new coordinates. We have (9), 
(34), 
and
$
\bar k_E(\bar x)=k_E(x)
\  ,
\ \
\bar F( \bar x,y)\ =\ \frac{\partial \phi}{\partial x}(x)\,  
F(x,y)
$.
This implies
\begin{eqnarray*}
\bar F(\hat {\bar {x}},y)&=&
\frac{\partial \phi}{\partial x}(\xhat)\left[  
f(\xhat) \;-\; k_E(\xhat)\,  P(\xhat)^{-1}
\frac{\partial h}{\partial x}(\xhat) ^\top \times \right. \\
& & \left.
\frac{\partial \delta }{\partial y_1}(\h(\xhat),y)
\right]
\  ,
\\
&=&\bar f(\bar {\hat {x}})
\;-\; \bar k_E(\bar x)
\left(\frac{\partial \phi}{\partial x}(\xhat) P(\xhat)^{-1}
\frac{\partial \phi}{\partial x}(\xhat)^\top\right)\times\\
& &\left[\frac{\partial \phi}{\partial x}(\xhat)^\top\right]^{-1}
\frac{\partial h}{\partial x}(\xhat) ^\top
\frac{\partial \delta }{\partial y_1}(\h(\xhat),y)
\  ,
\\
&=&
\bar f(\bar {\hat {x}})
\;-\; \bar k_E(\bar x)
\bar P(\bar {\hat {x}})^{-1}
\frac{\partial \bar h}{\partial \bar x}(\bar {\hat {x}})^\top
\frac{\partial \delta }{\partial y_1}(\bar h(\bar {\hat {x}},y))
\  .
\end{eqnarray*}
Therefore, the expression of the observer remains the same after the change 
of coordinates.
\end{IEEEproof}

\section{Conclusion}
If for a Riemannian metric $P$ and an observer such that 
the distance between estimated state and system state
decreases along the solutions, then the Lie derivative of $P$ along 
the systems solutions satisfies the inequality
in Theorem~\ref{thm3} involving the output function.
Also, the satisfaction of
such an inequality together with the existence of upper and 
lower bounds for
$P$
(see (\ref{eqn:Pbounds})) imply
detectability of the linear time-varying systems
obtained from linearizing the given system \eqref{eqn:Plant1} along its solutions.
Moreover, we have seen how the geodesic convexity of the output function level sets
is
necessary if the observer has an infinite gain margin
and, in a general situation, when the
Riemannian distance between estimated state and system state
decreases along the solutions of (\ref{eqn:Interconnection}).

Conversely, from the data of a
Riemannian metric
satisfying 
the necessary conditions in
Theorem~\ref{thm3}
and (\ref{eqn:Pbounds}),
and when the level sets of the output function are
totally geodesic,
we showed how to construct,
for the single output case,
an observer guaranteeing convergence of 
the estimated state to the system state, 
semiglobally with respect to zero estimation error set $\A$.

Also, although in Section \ref{sec1}
we have given an 
expression of an observer, at this 
time, we consider this only as an existence result and not as an 
observer design interesting for application. 
Actually we have
investigated mainly only the possibility and interest of studying observer 
convergence via a Riemannian metric, crystallizing the idea of using a  
contraction property. 
In a companion paper,
we focus on observer design,
where we study several scenarios in which
it is possible to
construct a Riemannian metric
 satisfying the desired inequality on its Lie derivative and making 
 the level sets of the output function possibly totally geodesic.

As a final remark, we observe that extensions of the results to 
nonautonomous systems, in particular those with inputs, seem possible using the proof techniques proposed here.
Also time scaling exploiting the concept of unbounded observability, 
as in \cite{Astolfi.Praly.06}, is expected  to be useful in relaxing 
the system completeness assumption.

\appendix
\label{sec:appendix}

\subsection{A necessary condition for completeness}
The following lemma provides conditions on $P$ that
guarantee that geodesics can be maximally extended to $\reals$.
\begin{lemma}
\label{lem1}
{\it
Suppose that a 
symmetric
covariant two-tensor
 $P:\reals^n\to \reals^{n\times n}$ satisfies
\begin{equation}
\label{3}
0\: <\:  P(x)
\quad \forall x\in \reals ^n\  ,\quad 
\lim_{r \to \infty }r^2 \plower(r) \:=\: +\infty,
\end{equation}
where, for any positive real number $r$,
$
\plower(r)\;=\; \min_{x:|x| \leq r}
\lambda _{\min}\left(P(x)\right).
$
Then, with $P$ as Riemannian metric
 on $\reals^n$,
any geodesic can be maximally 
extended to $\reals$.
}
\end{lemma}
\begin{IEEEproof}
Let $x_1$ and $x_2$ be any point in the ball $B_r$ in $\reals^n$ centered 
at the origin and with radius $r$. The Euclidean distance $|x_1-x_2|$ satisfies
$
\int_{s_1}^{s_2} \left|\frac{d  \gamma }{ds}(s )\right|  ds
\; \geq \;   |x_1-x_2|
$,
where $\gamma $ is any
piecewise 
$C^1$ path
between $x_1$ and $x_2$.
Using \eqref{eqn:Ldefn},
this implies
that, for any positive number $r$, 
\begin{equation}
\label{LP3}
\left.
\vrule height 1em depth 0.5em width 0pt
L(\gamma )\right|_{s_1}^{s_2}
\geq \, \sqrt{\plower(r)} 
\,  \int_{s_1}^{s_2} \left|\frac{d  \gamma }{ds}(s )\right|  ds
\, \geq \, \sqrt{\plower(r)} \:  |x_1-x_2|
\,  .
\end{equation}

Let $\gamma $ be any
normalized
geodesic maximally defined on $(\sigma 
_-,\sigma _+)$. By definition, it satisfies
\begin{equation}
\label{2}
\frac{d\gamma }{ds}(s) ^\top P(\gamma (s)) \,  \frac{d\gamma }{ds}(s)\;=\; 1
\qquad \forall s \in (\sigma _-,\sigma _+)
\  .
\end{equation}
Let $[s_1,s_2]$ be any closed interval contained in 
$(\sigma _-,\sigma _+)$. The function $\gamma :[s_1,s_2]\to \reals ^n$ 
is bounded (with the Euclidean norm).
We denote
$
r_{[s_1,s_2]}\;=\; \max_{s\in [s_1,s_2]}|\gamma (s)|.
$
By continuity, there exists $s_{12}$ in $[s_1,s_2]$ 
satisfying
$
r_{[s_1,s_2]}\;=\; |\gamma (s_{12})|.
$
Then, from (\ref{LP3}) and (\ref{2}), we obtain
\begin{equation}
\label{1}
\sqrt{\plower(|\gamma (s_{12})|)}\:  |\gamma (s_{12})-\gamma (s_2)|
\; \leq \;\left.\vrule height 1em depth 0.5em width 0pt
L(\gamma )\right|_{s_{12}}^{s_2}\;=\; |s_{12}-s_2|
\  .
\end{equation}

Because $(\sigma _-,\sigma _+)$ is the maximal interval of 
definition of $\gamma $, if $\sigma _-$ is finite, we must have

\noindent
$
\lim_{s_1\to \sigma _-} 
\left|\left(\gamma (s_1),\frac{d \gamma }{ds}(s_1)\right)\right|\;=\; +\infty 
$.
Now in the case where we have
$
\lim_{s_1\to \sigma _-} |\gamma (s_1)|\;=\; +\infty 
$
the definition of $s_{12}$ implies
$
\lim_{s_1\to \sigma _-}
\max_{s\in [s_1,s_2]}|\gamma (s)| =
\lim_{s_1\to \sigma _-} |\gamma (s_{12})|\;=\; +\infty.
$
Then, with assumption (\ref{3})
and (\ref{1}),
we get
\begin{eqnarray*}
|\sigma _--s_2|\:  & \geq & \:
\lim_{s_1\to \sigma _-}
\sqrt{\plower(|\gamma (s_{12})|)}\:  |\gamma (s_{12})-\gamma (s_2)|\:\geq \: \\
& & \hspace{-0.4in} 
\lim_{s_1\to \sigma _-}
\sqrt{\plower(|\gamma (s_{12})|)}(
|\gamma (s_{12})|-|\gamma (s_2)|
)
 \:\geq \:+\infty 
\  .
\end{eqnarray*}
This
is a contradiction.
Then, we are left with the case
$
\lim_{s_1\to \sigma _-} 
\left|\frac{d \gamma }{ds}(s_1)\right|\;=\; +\infty$.
But this contradicts (\ref{2}) since
we just established that $\gamma $ is 
bounded on $(\sigma _-,s_2)$, which, with \eqref{3},
implies that $P \circ \gamma $ is bounded away from $0$.

The same arguments apply to show that $\sigma _+ = +\infty$.
\end{IEEEproof}
\par\vspace{0.8412em}
\subsection{
On totally geodesic sets and property H3'
}

\begin{proposition}
{\it
Let $P$ be a complete Riemannian metric
on $\reals^n$ and $\mathcal{C}$ be a geodesically 
convex
subset of 
$\reals^n$.
\begin{enumerate}
\item
\label{itemzero}
If there exists $x_0$ in $\mathcal{C}$ satisfying
$
\frac{\partial h}{\partial x}(x_0)\;=\;  0
$
and all the sets $\mathfrak{H}(y)\cap \mathcal{C} $ for $y$ in 
$h(\mathcal{C})$
are totally geodesic then $h$ is constant on $\mathcal{C}$.
\item
\label{itemone}
Let $\mathcal{O}$ be the following open subset of $\reals^n$:
\begin{equation}\label{eqn:RankCondition}
\mathcal{O}\;=\; \left\{x\in\mbox{\rm int}(\mathcal{C})\,  :\: \mbox{\rm Rank}\left(\frac{\partial 
h}{\partial x}(x)
\right)=m
\right\}
\  .
\end{equation}
If all the sets  $\mathfrak{H}(y)\cap \mathcal{C}$ for $y$ in $h(\mathcal{C})$
are totally geodesic then we have,
for all $(j,k,l)$ and all $x \in \mathcal{O}$,
\begin{equation}
\label{LP55}
\begin{array}{lll} \displaystyle
\frac{\partial ^2h_j}{\partial x_k \partial x_l}(x)
-
\sum_{i=1}^n
\frac{\partial h_j}{\partial x_i}(x)\,  \Gamma _{kl}^i(x)
\displaystyle 
& = & \\
& & \hspace{-1.8in} \displaystyle
\sum_{i=1}^n 
\left(
g_{jik}(x)
\frac{\partial h_i}{\partial x_l}(x)
+
g_{jil}(x)
\frac{\partial h_i}{\partial x_k}(x)
\right)
\  ,
\end{array}
\end{equation}
where 
$
g_{jik}
:\mathcal{O}\to \reals$ are continuous
arbitrary functions
and $\Gamma_{kl}^i$
are the Christoffel symbols
\begin{equation}\label{eqn:ChristoffelSymbols}
\begin{array}{lll} \displaystyle
\Gamma^i_{kl}(x)
 & =&  \displaystyle
\frac{1}{2} \sum_{m = 1}^n
\left(P(x)^{-1}\right)_{im} \displaystyle
\left(
\frac{\partial P_{mk}}{\partial x_l}(x) \right.
\\
& & \hspace{0.7in} \displaystyle
\left.
+
\frac{\partial P_{ml}}{\partial x_k}(x)
-
\frac{\partial P_{kl}}{\partial x_m}(x)
\right).
\end{array}
\end{equation}

Conversely, if (\ref{LP55}) holds for any $x$ 
in $\mathcal{C}$, then
all the sets  $\mathfrak{H}(y)\cap \mathcal{C}$ for $y$ in 
$h(\mathcal{C})$
are totally geodesic.
\end{enumerate}
}
\end{proposition}

{\it Proof of item \ref{itemzero}:} The set $\mathcal{C}$ being geodesically convex, for any $x$ 
there exists a minimal geodesic $\gamma ^*$
between $x_0=\gamma ^*(0) $ and $x=\gamma ^*(s) $
satisfying
$
\gamma ^*(\sigma )\in\mathcal{C}\  \forall \sigma \in [0,s]$. Since we have
$
\frac{\partial h}{\partial 
x}(x_0)\frac{d\gamma ^*}{ds}(0)= 0
$
and the set $\mathfrak{H}(h(x_0))\cap\mathcal{C}$ is totally geodesic, we get
$
h(x)= h(x_0),
$
$x$ being arbitrary in $\mathcal{C}$, $h$ must be constant on 
$\mathcal{C}$.

{\it Proof of item \ref{itemone} (necessity):}
If $\mathcal{O}$ is empty, the statement holds vacuously.
If $\mathcal{O}$ is nonempty,
let $x$ be in $\mathcal{O}$. It is in the 
totally geodesic set
$\mathfrak{H}(h(x))\cap \mathcal{C}$. Then, for any $v$ in
$\reals ^n$ satisfying
\begin{equation}
\label{LP52}
\frac{\partial h}{\partial x}(x)\,  v\;=\; 0
\quad ,\qquad
v^\top P(x)\,  v\;=\; 1
\  ,
\end{equation}
consider a geodesic $\gamma $ satisfying
\begin{equation}\label{eqn:AuxGeoConvexityProof}
\gamma (0)\;=\; x\quad ,\qquad 
\frac{d\gamma }{ds}(0)\;=\; v
\end{equation}
with values in $\mathcal{C}$ on an interval
$(\sigma _-,\sigma _+)$.
We have
$
h(\gamma (s)) \;=\; 0$ for all $s \in (\sigma _-,\sigma _+)$.
This implies that we have
\begin{equation}\label{eqn:GeoConvexityCondition}
\frac{d h \circ \gamma }{ds}(0)\;=\; 
\frac{d^2 h \circ \gamma }{ds^2}(0)\;=\; 0
\  .
\end{equation}
But, with the geodesic equation, if we let
$
Q_{jkl}(x)\;=\; \frac{\partial ^2h_j}{\partial x_k \partial x_l}(x)
\;-\;
\sum_{i=1}^n
\frac{\partial h_j}{\partial x_i}(x)\,  \Gamma _{kl}^i(x)
$,
we have
\begin{equation}
\label{LP56}
\frac{d^2 h_j \circ \gamma }{ds^2}(s)\;=\; 
\sum_{k=1}^n\sum_{l=1}^n
Q_{jkl}
(\gamma (s))\frac{d \gamma _k}{ds}(s) 
\frac{d \gamma _l}{ds}(s) 
\  .
\end{equation}
Then, using
\eqref{eqn:AuxGeoConvexityProof}
and
\eqref{eqn:GeoConvexityCondition}, we have
\begin{equation}
\label{LP53}
\sum_{k=1}^n\sum_{l=1}^nQ_{jkl}(x)v_kv_l\;=\; 0
\qquad \forall j  \in \{1,2,\ldots,m\},
\end{equation}
where $v_k$ is the $k$th component of $v$.
Hence, we have established
$\sum_{k=1}^n\sum_{l=1}^nQ_{jkl}(x)v_kv_l\;=\; 0$
for all $(j,v=(v_k),x)\,  :\: 
j  \in \{1,2,\ldots,m\}
\,  ,\: \frac{\partial h}{\partial x}(x)v=0
\,  ,\: 
x\in \mathcal{O}
\  .$
The result follows from the S-Lemma (see \cite{Polik-Terlaky.07} for 
instance). In particular, we can pick
the functions
$g_{jik}(x)$
satisfying 
\eqref{LP55}
as,
for each $j$,
the entries of the matrix
\begin{eqnarray*}
\begin{array}{l}\displaystyle
\left[\frac{\partial h}{\partial x}(x)
\frac{\partial h}{\partial x}(x)^\top
\right]^{-1}\frac{\partial h}{\partial x}(x)
Q_{j\bullet\bullet}(x)
\times \\ \displaystyle
\hspace{0.5in}\left(I-
\frac{\frac{\partial h}{\partial x}(x)^\top[\frac{\partial h}{\partial x}(x)\frac{\partial h}{\partial x}(x)^\top]^{-1}\frac{\partial h}{\partial x}(x)}{2}\right)
\  .
\end{array}
\end{eqnarray*}

{\it Proof of item \ref{itemone} (sufficiency):}
For any $y$ in $h(\mathcal{C})$, let $(x,v)$ be any pair in 
$(\mathfrak{H}(y)\cap \mathcal{C})\times\reals^n$
satisfying
$
h(x) =  y$,
$\frac{\partial h}{\partial x}(x)\,  v\;=\; 0$,
$v^\top P(x)\,  v\;=\; 1$
and let $\gamma $ be any geodesic satisfying
$
\gamma (0)\;=\; x,\ \frac{d\gamma }{ds}(0)\;=\; v.
$
Let $J_\gamma $ be
the maximal interval containing $0$ so that $\gamma (J_\gamma )$ 
is contained in $\mathcal{C}$.
If $J_\gamma $ is reduced to a point, there is nothing 
to prove. If not $J_\gamma $ is an interval
with a non empty interior.
Then, with (\ref{LP56}) and (\ref{LP55}), for any
interior point  $s$  of $J_\gamma $, we have, 
for each $j$ in $\{1,\ldots,m\}$,
\begin{eqnarray*}
\frac{d}{ds}\frac{dh_j\circ \gamma }{ds}(s) &=&
\sum_{k=1}^n\sum_{l=1}^n
Q_{jkl}
(\gamma (s))\frac{d \gamma _k}{ds}(s) 
\frac{d \gamma _l}{ds}(s)  \\
&=&
2
\sum_{i=1}^n
\left[\sum_{k=1}^n
g_{jik}
(\gamma (s))\frac{d \gamma _k}{ds}(s) \right]
\frac{d h_i\circ \gamma }{ds}(s)
\  .
\end{eqnarray*}
Let $M$ be the matrix with entries $M_{ji}$ defined as,
$
M_{ji}(s)
= 2\left[\sum_{k=1}^n g_{jik}(\gamma (s))\frac{d \gamma _k}{ds}(s) 
\right],
$
for each $s\in \textrm{int}(J_\gamma )$.
The linear time varying system
$
\frac{dz}{ds}=M(s) z
$
has unique solutions. The only one satisfying $z(0)=0$ is 
identically $0$. So with the uniqueness of the solution of the 
geodesic equation we must also have
$\; 
\frac{d h_j\circ \gamma }{ds}(s)=0
\   \forall s\in \textrm{int}(J_\gamma )\; 
$
and therefore
$
h_j(\gamma (s)) = y_j$
for each $s\in \textrm{int}(J_\gamma )$ and each $j$.
Also, by continuity, if the upper bound $\sigma _+$ (respectively 
lower bound $\sigma _-$) of $J_\gamma $ is in $J_\gamma $,
then
we have also
$
h_j(\sigma _+)\;=\; y_j\ (\textrm{respectively}\quad h_j(\sigma _+)\;=\; y_j\; ).
$

\begin{proposition}
\label{prop2}
{\it
Let $P$ be a complete Riemannian metric
on $\reals^n$ and $\mathcal{C}$ be a geodesically 
convex
subset of 
$\reals^n$.
\begin{enumerate}
\item
If property H3' holds then
all the sets  $\mathfrak{H}(y)\cap 
\mathcal{C}$
 for $y$ in $h(\mathcal{C})$ are
\begin{enumerate}
\item
\label{itemfive}
 totally geodesic,
\item
\label{itemsix}
and geodesically convex.
\end{enumerate}
\item
If $m=1$
and all the sets  $\mathfrak{H}(y)\cap 
\mathcal{C}$
 for $y$ in $h(\mathcal{C})$ are totally geodesic then 
\begin{enumerate}
\item
\label{itemtwo}
they are all geodesically convex,
\item
\label{itemthree}
and property H3' holds with
$$
\delta (y_1,y_2)\;=\; |y_1-y_2|^2.
$$
\end{enumerate}
\end{enumerate}
}
\end{proposition}

{\it Proof of item \ref{itemfive}:}
Let $(x,v)$ be an arbitrary pair  in $\mathcal{C}\times \reals^n$
satisfying
\begin{equation}
\label{LP61}
\frac{\partial h}{\partial x}(x)\,  v\;=\; 0
\quad ,\qquad 
v^\top P(x)\,  v\;=\; 1
\  .
\end{equation}
Consider the geodesic
$\gamma_v$ satisfying
\begin{equation}
\label{LP62}
\gamma _v(0)\;=\; x\quad ,\qquad \frac{d\gamma _v}{ds}(0)\;=\; v
\  .
\end{equation}
Since $P$ is complete,
 ${\gamma_v}$ is defined on $(-\infty ,+\infty )$.
Let $J_{\gamma_v}$ be the maximal interval containing $0$ so that $\gamma_v (J_{\gamma_v})$ 
is contained in $\mathcal{C}$.

If $J_{\gamma_v}$ is reduced to a point, there is nothing to prove.
In the other case, for the sake of getting a contradiction, assume that $h$ is not
constant along this geodesic on $J_{\gamma_v}$, i.e., there exists $s_0$ in $J_{\gamma_v}$, say positive, satisfying
$
h(\gamma _v( s_0))\neq h(x)$, $\gamma _v(\sigma )\in\mathcal{C}$ for all $\sigma \in [0, s_0].
$
Let $s_1$ be the infimum of the real numbers $s$ in $[0, s_0]$ 
satisfying
$
h(\gamma _v(s))\; \neq \;  h(x).
$
By continuity $s_1$ is in $[0, s_0)$ and  we have
$
h(\gamma _v(s_1))= h( x).
$
Also, the definition of $s_1$ implies that,
for any  $\varepsilon $ in $(0,s_0-s_1]$, there exits 
$s_\varepsilon $ in $[s_1,s_1+\varepsilon ]$ such that
$
h(\gamma _v(s_\varepsilon ))\neq  h(\gamma _v(s_1)).
$
Also, when $s_1\neq 0$, the function $s\mapsto h(\gamma _v(s))$ being
constant on $[0,s_1]$, we have
\begin{equation}
\label{LP60}
\frac{\partial h}{\partial x}(\gamma _v(s_1))\,  \frac{d\gamma 
_v}{ds}(s_1)\;=\; 0
\  .
\end{equation}
Note that, with (\ref{LP61}) and (\ref{LP62}), the same holds when $s_1=0$.

Now let $\mathcal{B}_\varepsilon (\gamma _v(s_1))$ be a geodesic ball 
centered at $\gamma _v(s_1)$ with geodesic radius $\varepsilon $ sufficiently small 
to ensure that each geodesic between $\gamma _v(s_1)$ and any point 
in this ball is minimal. See \cite[Theorem VI.7.2]{Boothby.75}. With 
$s_\varepsilon $
associated with $\varepsilon $ as shown
above, we define a function $\gamma ^*$ as
$
\gamma ^*(s)\;=\; \gamma _v(s_\varepsilon -s)$
for all $s\in [0,s_\varepsilon -s_1]$.
It is a minimal geodesic between $\gamma ^*(0)=\gamma _v(s_\varepsilon )$
and $\gamma ^*(s_\varepsilon -s_1)=\gamma _v(s_1)$ satisfying
$
\gamma ^*(s)\in\mathcal{C}\cap \mathcal{B}_\varepsilon (\gamma _v(s_1))$
for all $s \in [0,s_\varepsilon -s_1]$
and 
$
h(\gamma ^*(0))\; \neq \;h( \gamma ^*(s_\varepsilon -s_1)).
$
So, according to H3', we have
$$
\frac{d}{ds}\delta (h(\gamma ^*(s)),h(\gamma ^*(0)))\; >\; 0$$
for all $s \in (0,s_\varepsilon -s_1]$.
In particular, we have
\begin{eqnarray*}
\frac{\partial \delta }{\partial y_1}
(h(\gamma ^*(s_\varepsilon -s_1)),h(\gamma ^*(0)))\times & & \\
& & \, \hspace{-1in} 
\frac{\partial h}{\partial x}(\gamma ^*(s_\varepsilon -s_1))
\,  
\frac{d \gamma ^*}{ds}(s_\varepsilon -s_1) \; >\; 0.
\end{eqnarray*}
But (\ref{LP60}) leads 
to a contradiction since
\begin{eqnarray*}
\frac{\partial h}{\partial x}(\gamma ^*(s_\varepsilon -s_1))
\,  
\frac{d \gamma ^*}{ds}(s_\varepsilon -s_1) 
&= & -
\frac{\partial h}{\partial x}(\gamma _v(s_1))
\,  
\frac{d \gamma _v}{ds}(s_1)\\
&= & 0.
\end{eqnarray*}

{\it Proof of item \ref{itemsix}:}
Let $(x_1,x_2) \in \mathcal{C}\times\mathcal{C}$ be any arbitrary pair of points satisfying
$
h(x_1) = h(x_2)= y.
$
Since $\mathcal{C}$ is geodesically convex, there exists a 
minimal geodesic $\gamma ^*$
between $x_1=\gamma ^*(s_1)$ and $x_2=\gamma ^*(s_2)$
satisfying
$
\gamma ^*(s)\in\mathcal{C}$
for all $s\in [s_1,s_2]$.
We have
$
\delta (h(\gamma ^* (s_2)),h(\gamma ^* (s_1)))
=
 \int_{s_1}^{s_2}
\frac{\partial \delta }{\partial y_1}(h(\gamma ^*(s)),h(\gamma ^*(s_1)))$
$\,  \frac{d\,  h\circ \gamma ^* }{ds}(s)\,  ds
$.
But (\ref{LP39}) implies the left-hand side of this equation is zero if and only if
we have
$
h(\gamma ^* (s))= h(\gamma ^* (s_1))$
for all $s\in [s_1,s_2]$,
that is, the geodesic $\gamma ^* $ remains in the set
$\mathfrak{H}(h(x_1))\cap \mathcal{C}$ for all $s$ in $[s_1,s_2]$.

{\it Proof of item \ref{itemtwo}:}
Let $(x_1,x_2) \in \mathcal{C}\times\mathcal{C}$ be any arbitrary pair of points satisfying
$h(x_1) =h(x_2)
=y
$.
Since $\mathcal{C}$ is geodesically convex, there exists a 
minimal geodesic $\gamma ^*$
between $x_1=\gamma ^*(s_1)$ and $x_2=\gamma ^*(s_2)$
satisfying
$
\gamma ^*(s)\in\mathcal{C}$ for all $s\in [s_1,s_2]
$.  
For the sake of getting a contradiction, assume 
 that $\mathfrak{H}(y)\cap \mathcal{C}$
is not geodesically convex.
Then,
there exists $\hat s \in [s_1,s_2]$
such that $\gamma ^*(s)
\not\in \mathfrak{H}(y)\cap \mathcal{C}
$. But $\gamma^*(\hat s)$ being in $\mathcal{C}$, this implies
$
|h(\gamma ^*(\hat s))- h(x_1)|^2 \neq  0.
$
By continuity and compactness, the function $s\in [s_1,s_2]\mapsto 
|h(\gamma ^*(s))-h(x_1)|^2
$ admits a maximum at some $s_{\max}$ in $(s_1,s_2)$ 
and, hence
\begin{eqnarray}
\label{LP36}
h(\gamma ^*(s_{\max}))
\neq 
 \! h(x_1),
\end{eqnarray}
\begin{eqnarray}
\nonumber 
\begin{array}{lll}
\displaystyle
\left(h(\gamma ^*(s_{\max}))-h(x_1)\right)^\top
\frac{dh\circ\gamma ^*}{ds}(s_{\max})
&=& \\  \nonumber
\\ \nonumber
& &  \hspace{-2.5in}
\displaystyle
\left(h(\gamma ^*(s_{\max}))-h(x_1)\right)^\top
\frac{\partial h}{\partial x}(\gamma ^*(s_{\max}))
\frac{d\gamma ^*}{ds}(s_{\max})=0\  .
\end{array}
\end{eqnarray}
When the dimension $m$ of outputs is one, this implies
$
\frac{\partial h}{\partial x}(\gamma ^*(s_{\max}))
\frac{d\gamma ^*}{ds}(s_{\max})= 0.
$
Since the set $\mathfrak{H}(h(\gamma ^*(s_{max})))\cap\mathcal{C}$
is totally geodesic and $\gamma ^*$ takes its values in $\mathcal{C}$ on 
the interval $[s_1,s_2]$ containing $s_{max}$, we conclude that
$\gamma ^*$ takes actually its values in
$\mathfrak{H}(h(\gamma ^*(s_{\max})))\cap\mathcal{C}$
on $[s_1,s_2]$. This contradicts (\ref{LP36}), 
and so  $\mathfrak{H}(y)\cap \mathcal{C}$
must be geodesically convex.

{\it Proof of item \ref{itemthree}:}
Let $(\xhat ,x)$ be an arbitrary pair of points
in $\mathcal{C}\times\mathcal{C}$ satisfying
$
h(\xhat)\; \neq \; h(x)
\  .$
Since $\mathcal{C}$ is geodesically convex, there exists a 
minimal geodesic $\gamma ^*$
between $x=\gamma ^*(0)$ and $\xhat =\gamma ^*(\hat s)$
satisfying
$\gamma ^* (s)\in\mathcal{C}$ for all $s\in [0,\hat s]$.
Assume there exists $s$ in $[0,\hat s]$ satisfying
$\frac{d\,  h\circ\gamma ^* }{ds}(s)=
\frac{\partial h}{\partial x}(\gamma ^* (s))\,  \frac{d\gamma ^* }{ds}(s)= 0$,
$\frac{d\gamma ^* }{ds}(s)^\top P(\gamma ^* (s))\,  \frac{d\gamma ^* }{ds}(s)= 1$.
Then,
since $\mathfrak{H}(h(\gamma ^* (s))\cap \mathcal{C}$ is totally geodesic,
and $\gamma ^* $ takes its values in $\mathcal{C}$ on $[0,\hat s]$,
we have
$
h \circ \gamma ^* (s) = h \circ \gamma ^* (0)= h(x)$ 
for all $s\in [0,\hat s]$
which contradicts
$
h \circ \gamma ^* (\hat s)\;=\; h(\xhat)\neq h(x)$.
Then, $\displaystyle \frac{d\,  h\circ\gamma ^* }{ds}$ has a constant sign. But,
since we have
$
h(\hat x)-h(x)= \int_0^{\hat s}
\frac{d\,  h\circ\gamma ^* }{ds}(s) \,  ds$,
this sign must be the same as the one of $h(\hat x)-h(x)$.
We conclude that we have
$$
\frac{d}{ds}|h(\gamma (s))-h(\gamma (0))|^2
= 
[h(\gamma (s))-h(\gamma (0))]\,  \frac{d\,  h\circ\gamma ^* }{ds}(s)
 > 0$$
for all $s \in (0,\hat s]$.

\balance
\bibliographystyle{plain} 
\bibliography{ObserverGeodesic}

\begin{thebibliography}{10}

\bibitem{Aghannan.Rouchon.03}
N.~Aghannan and P.~Rouchon.
\newblock An intrinsic observer for a class of {L}agrangian systems.
\newblock {\em IEEE Trans. Automatic Control}, 48(6):936--945, 2003.

\bibitem{Astolfi.Praly.06}
A.~Astolfi and L.~Praly.
\newblock Global complete observability and output-to-state stability imply the
  existence of a globally convergent observer.
\newblock {\em Mathematics of Control, Signals, and Systems {(MCSS)}},
  18(1):32--65, 2006.

\bibitem{Bonnabel.07}
S.~Bonnabel.
\newblock {\em Observateurs Asymptotiques Invariants}.
\newblock PhD thesis, Sp\'{e}\-cia\-li\-t\'{e} ñMa\-th\'{e}\-ma\-ti\-ques et
  Automatiqueî. \'Ecole des Mines de Paris, 2007.

\bibitem{Bonnabel.10}
S.~Bonnabel.
\newblock A simple intrinsic reduced-observer for geodesic flow.
\newblock {\em IEEE Transactions on Automatic Control}, 55(9):2186--2191, 2010.

\bibitem{Boothby.75}
W.~M. Boothby.
\newblock {\em An Introduction to Differentiable Manifolds and {R}iemannian
  Geometry}.
\newblock Academic Press, 1975.

\bibitem{Cartan.51}
E.~Cartan.
\newblock {\em Geometry of {R}iemannian spaces}.
\newblock Maths Sci Press, second edition, 1951.

\bibitem{Demidovich.61-62}
B.~P. Demidovich.
\newblock Dissipativity of a nonlinear system of differential equations.
\newblock {\em Ser. Mat; Mekh. Part I.6 (1961); Part II.1, 3-8(1962) (in
  Russian).}, pages 19--27.

\bibitem{DoCarmo.92}
M.P. DoCarmo.
\newblock {\em Riemannian Geometry}.
\newblock Birkhauser, Boston, 1992.

\bibitem{Filippov.88}
F.~Filippov.
\newblock {\em Differential equations with discontinuous right hand sides}.
\newblock Kluwer Academic Publishers. Mathematics and Its Applications, 1988.

\bibitem{Hartmann.64}
P.~Hartmann.
\newblock {\em Ordinary differential equations}.
\newblock Wiley, 1964.

\bibitem{Hicks.65}
N.~J. Hicks.
\newblock {\em Notes on Differential Geometry}.
\newblock Van Nostrand Publishing Company, Princeton, New Jersey, 1965.

\bibitem{Polik-Terlaky.07}
T.~Terlaky I.~P\'{o}lik.
\newblock A survey of the s-lemma.
\newblock {\em SIAM Review}, 49(3):371--418, september 2007.

\bibitem{Isac.Nemeth.08}
G.~Isac and S.~Z. N\'emeth.
\newblock {\em Scalar and asymptotic scalar derivatives : theory and
  applications}.
\newblock Springer, 4th edition, 2008.

\bibitem{Jouffroy.05.CDC}
J.~Jouffroy.
\newblock Some ancestors of contraction analysis.
\newblock {\em Proc. IEEE Conf. Dec. Control}, pages 5450--5455, 2005.

\bibitem{Kaboyashi.Nomizu.96}
S.~Kaboyashi and K.~Nomizu.
\newblock {\em Foundations of Differential Geometry}, volume~1.
\newblock Wiley Classics Library Edition, 1996.

\bibitem{Krener.Isidori.83}
A.~Krener and A.~Isidori.
\newblock Linearization by output injection and nonlinear observers.
\newblock {\em Systems \& Control Letters}, 3:47--52, 1983.

\bibitem{Lebedev.Cloud.05}
L.~P. Lebedev and M.~J. Cloud.
\newblock {\em Tensor analysis}.
\newblock Norwood Mass., 2005.

\bibitem{Lewis.71}
D.~C. Lewis.
\newblock Metric properties of differential equations.
\newblock {\em American Journal of Mathematics}, 71:294--312, 1949.

\bibitem{Lohmiller.Slotine.98.Automatica}
W.~Lohmiller and J.-J. Slotine.
\newblock On contraction analysis for nonlinear systems.
\newblock {\em Automatica}, 34(6):671--682, 1998.

\bibitem{Nemeth.98}
S.~Z. N\'emeth.
\newblock {\em Geometric aspects of Minty-Browder monotonicity}.
\newblock PhD thesis, E \"{o}tv\"{o}s Lor\'{a}nd University. Budapest., 1998.

\bibitem{Praly.01.NOLCOS.Observers}
L~Praly.
\newblock On observers with state independent error {L}yapunov function.
\newblock In {\em Proceedings of the 5th IFAC Symposium "Nonlinear Control
  Systems" (NOLCOS'01)}, 2001.

\bibitem{Price.40}
G.~B. Price.
\newblock Definitions and properties of monotone functions.
\newblock {\em Bull. Amer. Math. Soc.}, 46:77--80, 1940.

\bibitem{Rapcsak.97}
T.~Rapcs\'{a}k.
\newblock {\em Smooth Nonlinear Optimization in $\reals^n$}.
\newblock Kluwer Academic, 1997.

\bibitem{Sakai.96}
T.~Sakai.
\newblock {\em Riemannian geometry}, volume 149.
\newblock Translation of Mathematical monographs, American Mathematical Soc.,
  1996.

\bibitem{Spivak.79}
M.~Spivak.
\newblock {\em (A Comprehensive Introduction to) Differential geometry},
  volume~2.
\newblock Publish or Perish, Inc., 1979.
\newblock 2nd Edition.

\bibitem{Tsinias.90}
J.~Tsinias.
\newblock Further results on the observer design problem.
\newblock {\em Syst. Contr. Lett.}, 14:411--418, 1990.

\end{thebibliography}

\end{document}